\documentclass[12pt,reqno]{amsart}

\usepackage{amssymb,mathtools}

\usepackage{graphicx}
\usepackage[dvipsnames]{xcolor}
\usepackage[british]{babel}
\usepackage[T1]{fontenc}

\usepackage[margin=1in]{geometry}
\usepackage[autostyle]{csquotes}
\usepackage[style=alphabetic,url=true]{biblatex}
\usepackage{enumitem,float,tikz,listings,array,booktabs,longtable}

\usetikzlibrary{shapes,arrows,chains,intersections,calc,patterns}

\usepackage[colorlinks=true,
            linkcolor=NavyBlue,
            citecolor=NavyBlue,
            filecolor=NavyBlue,
            urlcolor=NavyBlue]{hyperref}

\theoremstyle{plain}
\newtheorem{theorem}{Theorem}
\newtheorem{corollary}[theorem]{Corollary}
\newtheorem{lemma}[theorem]{Lemma}
\newtheorem{proposition}[theorem]{Proposition}

\theoremstyle{definition}
\newtheorem{definition}[theorem]{Definition}
\newtheorem{example}[theorem]{Example}

\theoremstyle{remark}
\newtheorem{remark}[theorem]{Remark}

\newcommand{\defin}[1]{%
	\relax\ifmmode%
	\textcolor{blue}{#1}%
	\else\textcolor{blue}{\emph{#1}}%
	\fi%
}

\newcommand{\seqnum}[1]{\href{https://oeis.org/#1}{\underline{#1}}}

\DeclareMathOperator{\mb}{mb}
\DeclareMathOperator{\NMB}{Nmb}
\DeclareMathOperator{\nmb}{nmb}
\DeclareMathOperator{\MB}{Mb}
\DeclareMathOperator{\sueq}{suc}
\DeclareMathOperator{\Succeq}{Suc}
\DeclareMathOperator{\blocks}{bl}
\DeclareMathOperator{\maj}{maj}
\DeclareMathOperator{\des}{des}
\DeclareMathOperator{\Des}{Des}
\DeclareMathOperator{\Pk}{Pk}
\DeclareMathOperator{\Val}{Val}
\DeclareMathOperator{\Dasc}{Dasc}
\DeclareMathOperator{\Ddesc}{Ddesc}
\DeclareMathOperator{\pk}{pk}
\DeclareMathOperator{\val}{val}
\DeclareMathOperator{\dasc}{dasc}
\DeclareMathOperator{\ddesc}{ddesc}
\DeclareMathOperator{\Orb}{Orb}
\DeclareMathOperator{\Fix}{Fix}
\DeclareMathOperator{\Swap}{Swap}

\newcommand{\typeBnzb}[1]{{\Pi^0_{#1}}} 

\newcommand{\interl}{\ll}

\ExplSyntaxOn
\seq_new:N   \l__tb_blocks_seq
\clist_new:N \l__tb_elems_clist
\tl_new:N    \l__tb_tmp_tl
\int_new:N   \l__tb_val_int
\cs_new_protected:Npn \__tb_format_atom:n #1
{
	\tl_set:Nn \l__tb_tmp_tl {#1}
	\tl_trim_spaces:N \l__tb_tmp_tl
	\tl_if_blank:nTF { \l__tb_tmp_tl } { } 
	{
		\int_set:Nn \l__tb_val_int { \l__tb_tmp_tl }
		\int_compare:nNnTF { \l__tb_val_int } < { 0 }
		{ \overline{ \int_eval:n { -\l__tb_val_int } } }
		{ \int_eval:n { \l__tb_val_int } }
	}
}
\cs_new_protected:Npn \__tb_format_block:n #1
{
	\clist_set:Nn \l__tb_elems_clist {#1}
	\clist_map_inline:Nn \l__tb_elems_clist { \__tb_format_atom:n{##1}\, }
}
\NewDocumentCommand{\typeBPartition}{m}
{
	\ensuremath{
		\seq_set_split:Nnn \l__tb_blocks_seq {|} {#1}
		\seq_map_indexed_inline:Nn \l__tb_blocks_seq
		{
			\__tb_format_block:n{##2}
			\int_compare:nNnT {##1} < { \seq_count:N \l__tb_blocks_seq } {{\boldsymbol{\mid}\,}}
		}
	}
}
\ExplSyntaxOff

\title{Some Families of Type $B$ Set Partitions Counted by the Dowling Numbers}
\date{}

\author{Per Alexandersson}
\address{Department of Mathematics, Stockholm University, SE-106 91 Stockholm, Sweden}
\email{per.w.alexandersson@gmail.com}

\author{Fufa Beyene}
\address{Department of Mathematics, Kotebe University of Education, 31248 Addis Ababa, Ethiopia}
\email{fufa.beyene@aau.edu.et}

\author{Roberto Mantaci}
\address{IRIF, Université de Paris, 8 Place Aurélie Nemours, F-75013 Paris, France}
\email{mantaci@irif.fr}

\subjclass[2020]{Primary 05A05; Secondary 05A15, 05A19.}
\keywords{Type $B$ set partition, merging-free set partition, separated set partition, flattened Stirling permutations, flattened signed permutation, descent, $\gamma$-positivity, homomesic.}

\begin{document}

\begin{abstract}
    In this paper, we study type $B$ set partitions without zero block.
    Certain classes of these partitions, such as merging-free and separated 
    partitions (enumerated by the Dowling numbers), are investigated.
    We show that these classes are in bijection with type $B$ set partitions.
    The intersection of these two classes is also studied, and we prove that their block-generating polynomials
    are real-rooted.
    
    Finally, we study the descent statistics on the class of permutations obtained by flattening type $B$ merging-free partitions. Using the valley-hopping action,
    we prove the \(\gamma\)-positivity of the descent
    distribution and provide a combinatorial interpretation of the \(\gamma\)-coefficients.
    We also show that the descent statistic is homomesic under valley-hopping.
\end{abstract}

\maketitle

\tableofcontents

	\section{Introduction and preliminaries}
	
	Let $n$ be a fixed positive integer and let $\defin{[n]} \coloneqq \{1, 2, \dotsc, n\}$. 
	A \defin{type $A$ set partition} of $[n]$ is a collection of pairwise disjoint non-empty
	subsets of $[n]$ such that their union forms the whole set $[n]$.
	The number of set partitions of $[n]$ having $k$ blocks is
	the \defin{Stirling number of the second kind}. 
	In addition to the statistics of the number of blocks, many different statistics have been studied over set partitions.
	The numbers of \emph{merging blocks}, \emph{successions}, and \emph{singletons} are just a few of them, (see \cite{Be-Ma,Ca,Ma-Mu}).
	
	The study of type $B$ set partitions began with the work of Montenegro~\cite{Mo}, 
	where he investigates the type $B$ analog of non-crossing set partitions from which he constructs the type $B$ associahedron.
	Using this construction, Reiner~\cite{Re} develops a definition for type $B$ set partitions.
	The type $B$ set partitions of the set \defin{$\langle n\rangle \coloneqq \{0, \pm1, \pm2, \pm3, \dotsc, \pm n\}$} 
	are counted by the Dowling number \cite[\seqnum{A007405}]{Sl} (see also \cite{Bucketal}).
	Type $B$ Stirling number of the second kind counts type $B$ set partitions according to the 
	number of pairs of non-zero blocks \cite{Sa-Sw}.
	In \cite{Me-Ra}, Mez\"o and Ramirez enumerated type $B$ set partitions without singletons.
	In \cite{Ch-Wa}, Chen and Wang proved that the joint distribution of the number of singleton pairs 
	and the number of adjacency pairs is symmetric over type $B$ set partitions without zero block, 
	in analogy with the result of Callan for type $A$ set partitions.
	
	For any set $S$, the function $\sigma: S\mapsto [n]$ corresponds to the word 
	$\sigma(1)\sigma(2)\cdots\sigma(n)$.
	In particular, a permutation is a word with distinct symbols.
	
	The study of the distribution of different statistics in the set of type $A$ permutations of $[n]$ 
	or the \defin{symmetric group} has a long history~\cite{Eh-St}.
	Any statistic equidistributed with \emph{inversion number}, or \emph{major index},
	on the symmetric group is called \defin{Mahonian}.
	The classical type $A$ Eulerian number $A(n,k)$ (\seqnum{A123125} in OEIS) counts 
	type $A$ permutations of $[n]$ according to several statistics.
	For instance, \emph{descents}, \emph{ascents} and \emph{excedances} are a few of them.
	The Eulerian numbers satisfy the symmetric relation $A(n,k)=A(n,n-k-1)$.
	The joint distribution of descent number and major index
	for the symmetric group is described by a well-known identity of Carlitz.
	
	The real-rootedness and unimodality of polynomials are major themes which have occupied mathematicians for the past few decades. 
	The property of \(\gamma\)-positivity implies symmetry and unimodality. 
	Recall that  a polynomial \(f(x)=\sum_{i}a_ix^i\in\mathbb{R}[x]\) is called \defin{\(\gamma\)-positive}, if
	\[
	f(x)=\sum_{i=0}^{\lfloor\frac{n}{2}\rfloor}\gamma_ix^i(1+x)^{n-2i}
	\]
	for some \(n \in \mathbb{N} \) and non-negative reals \(\gamma_0, \gamma_1, \ldots, \gamma_{\lfloor\frac{n}{2}\rfloor}\). 
	\defin{The Eulerian polynomials}, 
	\[
	A_n(t)=\sum_{k=0}^{n-1}A(n,k)t^k,
	\]
	are the classical examples of real-rooted and \(\gamma\)-positive polynomials in combinatorics. 
	Furthermore, the roots of \(A_j(t)\) interlace the roots of \(A_{j+1}(t)\) 
	for all \(j\) (see \cite{Athanasiadis2018,Al-Na,Wagner1992}). 
	
	The organization of the paper is as follows. 
	
	In Section~\ref{sec:setpartwithoutzero}, we study type \(B\) set partitions without zero block. 
	We show that the joint distribution of the numbers of merging blocks and successions is symmetric. 
	This generalizes the study of the joint distribution of these statistics in type \(A\) 
	set partitions studied in~\cite{Beyeneetal}. 
	We show that the generating polynomials of the distribution of blocks of type \(B\) merging-free partitions, 
	Proposition~\ref{prop:interlacingRoots}, and type \(B\) merging-free separated partitions, 
	Theorem~\ref{thmnombsuc}, are real-rooted. 
		
	In Section~\ref{sec:flatStirlingperm}, we present a bijection between type \(B\) 
	set partitions over \([n]\) and flattened Stirling permutations over \([n]\). 
	
	In Section~\ref{sec:flatsignedperm}, we study the set \(\mathcal{R}_n^B\) of signed permutations 
	obtained by flattening all type \(B\) merging-free partitions. 
	We study the descent distribution on this set, and prove, in Theorem~\ref{thm:gammapos}, that 
	the polynomials
	\[
	\sum_{\sigma\in\mathcal{R}_n^B}x^{\des(\sigma)}
	\]
	are \(\gamma\)-positive by using a group action called \emph{valley-hopping} and we 
	provide a combinatorial interpretation of the $\gamma$-coefficients. 
	Also, we show that the descent statistic is homomesic under valley-hopping.
	\begin{figure}[!ht]
		\centering
		\begin{tikzpicture}[scale=0.5,
			every node/.style={font=\small},
			setlabel/.style={font=\bfseries},
			line/.style={line width=0.9pt},
			]
			\draw[line,black]  (-8,-4) rectangle (8,6);
			\def\r{2.6}
			\coordinate (A) at (-2.0, 0.0);   
			\coordinate (C) at ( 1.0, 0.0);   
			
			\begin{scope}
				\fill[blue!12]   (A) circle (\r);
				\fill[orange!12] (C) circle (\r);
				\draw[line,black]   (A) circle (\r);
				\draw[line,black] (C) circle (\r);
			\end{scope}
			
			\begin{scope}
				\clip (A) circle (\r);
				\fill[pattern=north east lines, pattern color=black!60] (-6,-6) rectangle (6,6);
			\end{scope}
			
			\begin{scope}
				\clip (A) circle (\r);
				\clip (C) circle (\r);
				\fill[pattern=dots, pattern color=black!60] (-6,-6) rectangle (6,6);
			\end{scope}
			
			\node[setlabel]   at (-2,6.5) {Type $B$ set partitions};
			\node[setlabel]   at ($(A)+(-1.0,3.1)$) {Merging-free};
			\node[setlabel] at ($(C)+( 1.0,3.1)$) {Separated};
			
		\end{tikzpicture}
		\caption{Overview of the families of type $B$ set partitions studied in this paper. All set-partitions considered are without zero block. Merging-free set partitions are studied in Section~\ref{sec:Merging-free}. Separated set partitions are studied in Section~\ref{sec:Separated}. 
		}
	\end{figure}
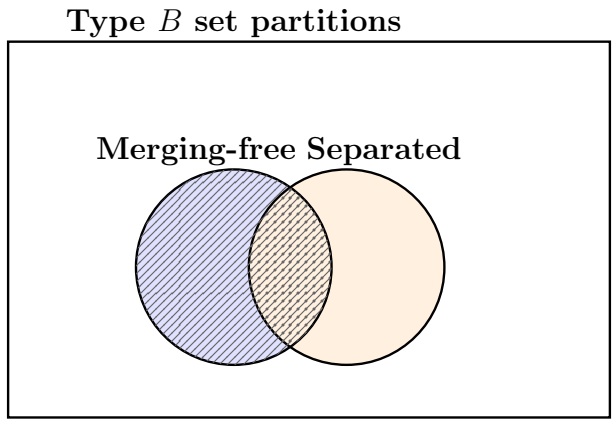


	\section{Background on type \texorpdfstring{\(B\)}{B} set partitions}
	
	\begin{definition}[Type $B$ set partition]
		
		A \defin{type $B$ set partition} $\pi=\pi_0 \mid \pi_1 \mid  \cdots \mid \pi_{2k}$ 
		is a set partition of $\langle n\rangle$ such that
		\begin{enumerate}
			\item for $i\ge1, \pi_{2i} = -\pi_{2i-1}$, where $-\beta \coloneqq \{-a: a\in\beta\}$;
			\item $\pi$ contains at most one \defin{zero block}, $\pi_0$,
			that has the property that $a\in \pi_0 \iff -a\in\pi_0$.
		\end{enumerate}
	\end{definition}
	Negative elements in a set partition are called \defin{barred} and we usually write $\overline{k}$ for $-k$.
	A block $\pi_i$ is called a \defin{singleton} if it contains exactly one element.
	\medskip
	
	Sagan~and~Swanson~\cite{Sa-Sw} defined the \defin{standard form} of type $B$ set partition as follows.
	Let $|\beta|=\{|a| : a\in\beta\}$, so that $|\pi_{2i}|=|-\pi_{2i-1}|$.
	For all $i$, we let $\defin{m_i} \coloneqq \min(|\pi_i|)$.
	Then $\pi$ is said to be in \defin{standard form} if
	\begin{enumerate}
		\item $m_{2i} \in \pi_{2i}$ for all $i$, and
		\item $0=m_0<m_2<m_4<\cdots<m_{2k}$.
	\end{enumerate}
	For instance, the standard form of the set partition
	\[
	\pi= \typeBPartition{-5 | -3, -6, 9 | -2, -7 | -1, 8 | 0, -4, 4 | 1, -8 | 2, 7 | 3, 6, -9 | 5 }
	\]
	is
	\[
	\pi= \typeBPartition{0, -4, 4 | 1, -8 | -1, 8 | 2, 7 | -2, -7 | 3, 6, -9 | -3, -6, 9| 5  | -5  }
	\]
	with $m_0=0$, $m_1=m_2=1$, $m_3=m_4=2$, $m_5=m_6=3$, $m_7=m_8=5$.
	
	Adler~\cite{Ad} denoted type $B$ set partition as follows: in the zero block,
	all negative elements are removed; in all blocks, the elements are sorted increasingly by
	absolute value; from any block pair we keep only the block with positive first element.
	
	\begin{example}
		Using Adler's notation, the above partition is written as
		\[
		\pi = \typeBPartition{0, 4 | 1, -8 | 2, 7 | 3, 6, -9 | 5}.
		\]
	\end{example}
	\medskip
	
	We mainly consider set partitions in their standard form written in Adler's notation.
	If $\pi=\pi_0 \mid \pi_1 \mid \dotsb \mid \pi_k$, then we let \defin{$\blocks(\pi)$}
	denote the number of blocks of $\pi$.
	Further, $m_i = \min(|\pi_i|)$, $0 \le i \le k$, is called the \defin{minimal element} of $\pi_i$
	and the element $M_i=\max(|\pi_i|)$ is called the \defin{maximal element} of $\pi_i$.
	
	\begin{example}
		Continuing with our example, $\blocks(\pi)=5$, the minimal elements are $m_0=0$, $m_1=1$, $m_2=2$, $m_3=3$, $m_4=5$,
		and the maximal elements are $M_0=4$, $M_1=8$, $M_2=7$, $M_3=9$, $M_4=5$.
	\end{example}
	
	Let \defin{$\Pi_n^B$} denote the set of all type $B$ set partitions over $\langle n\rangle$,
	and let
	\[
	\defin{\Pi_{n,k}^B} \coloneqq \{\pi\in\Pi_n^B: \blocks(\pi)=k\}.
	\]
	The cardinality of $\Pi_n^B$ is given by the \defin{Dowling number ($D_n$)}. 
	
	The number of type $B$ set partitions over $\langle n\rangle$ having $k$
	non-zero blocks satisfies the recurrence relation
	\begin{equation}\label{eq:typeBStirlingRecursion}
		S_B(n,k) = S_B(n-1,k-1)+(2k+1)S_B(n-1,k), \quad n\ge1,
	\end{equation}
	where $S_B(n,k)$ is the \defin{type $B$ Stirling numbers of the second kind}
	and $S_B(0,k)=\delta_{0,k}$ (Kronecker delta).


	\section{Type \texorpdfstring{$B$}{B} set partitions without zero block}\label{sec:setpartwithoutzero}
	
	Let \defin{$\Pi_n^0$} denote the set of type $B$ set partitions without zero block
	(or such that the zero block contains only $0$), and let \defin{$\Pi_{n,k}^0$} denote the set of type $B$
	set partitions without zero block and having $k$ blocks.
	Let \defin{$w_n$} denote the number of 
	type $B$ set partitions without zero block over $\langle n\rangle$.
	The first few terms of $w_n$ are
	\[
	1, \; 1, \; 3, \; 11, \; 49, \; 257, \; 1539, \; \dotsc.
	\]
	This is the sequence \cite[\seqnum{A004211}]{Sl}.
	\begin{proposition}
		For $n\ge0$, we have 
		\begin{equation}
			D_n=\sum_{j=0}^{n} \binom{n}{j} w_j.
		\end{equation}
	\end{proposition}
	\begin{proof}
		Suppose that the zero block of $\pi \in \Pi_n^B$ contains $j \in \{0,1,2,\dotsc,n\}$ elements.
		We can then choose the elements of the zero block of $\pi$ in $\binom{n}{j}$ ways.
		The remaining elements can be any type $B$ set partition of $n-j$ elements,
		and there are $w_{n-j}$ such partitions.
		Summing over all $j$, we get the identity above.
	\end{proof}
	
	\begin{corollary}
		The number $w_n$ of type $B$ set partitions without zero block is a
		binomial transform of the Dowling number $D_n$. That is,
		\begin{equation}
			w_n=\sum_{j=0}^{n}(-1)^{n-j}\binom{n}{j}D_j, n \geq 0.
		\end{equation}
	\end{corollary}
	\begin{proof}
		This follows from the basic properties of the binomial transform.
	\end{proof}

	\begin{proposition}
		For $n\ge0$, we have
		\begin{equation}
			w_n=\sum_{k=1}^{n}2^{n-k} S(n,k),
		\end{equation}
		where $S(n,k)$ is the Stirling number of the second kind with $S(n,0)=0$.
	\end{proposition}
	\begin{proof}
		There are $S(n,k)$ type $A$ set partitions over $[n]$ having $k$ blocks.
		From each of these partitions we obtain a type $B$ set partition over $\langle n\rangle$ by adding
		signs to the appropriate elements.
		Since the minimal elements in each of the $k$ blocks must remain positive,
		we have $n-k$ remaining elements that can be either positive or negative.
		There are $2^{n-k}$ ways to choose the signs for those elements.
		Thus, there are $2^{n-k} S(n,k)$ ways to construct a type $B$ set partition.
		By taking the sum over all possible values of $k$, we are done.
	\end{proof}

	The following proposition presents the analogous result
	of the recursion in \eqref{eq:typeBStirlingRecursion},
	but we now only count type $B$ set partitions without zero block.
	This is a special case of a more general result appearing in~\cite[Eq.~(1)]{Wang2014}.
	\begin{proposition}
		The number $w_{n,k}$ of type $B$ set partitions without
		zero block over $\langle n\rangle$ having $k$ blocks satisfies the recursion
		\begin{equation}
			w_{n,k}=2 k \cdot w_{n-1,k}+w_{n-1,k-1}, \text{ with }  w_{n,0}=\delta_{n,0}.
		\end{equation}
	\end{proposition}
	Wang also proves that the generating polynomials
	\begin{equation}
		T_n(x) = \sum_{k\geq 0} w_{n,k} x^k
	\end{equation}
	satisfy $T_n(x) = x T_{n-1}(x) + 2x T'_{n-1}(x)$
	and are real-rooted, see \cite[Thm. 2.1]{Wang2014}.

	\subsection{The distribution of the statistics of \texorpdfstring{$\mb$}{mb} and \texorpdfstring{$\sueq$}{suc}}\label{sec:mbSucc}
	
	The statistics $\mb$ and $\sueq$ over type $A$ 
	set partitions have been studied in \cite{Beyeneetal},
	and it has been shown that their joint distribution is symmetric.
	Now, we give a type $B$ analogue of the symmetric distribution of these statistics.

	\begin{definition}
		Let $\pi=\pi_1\mid \pi_2 \mid \dotsb \mid \pi_k \in \typeBnzb{n}$.
		
		For $2\leq i \leq k$, a block $\pi_i$ is said to be \defin{merging}
		if $\max(|\pi_{i-1}|)<\min(|\pi_i|)$.
		A partition with no merging block is called \defin{merging-free}.
		\medskip 
		
		The element $a \in \{2,3,\dotsc,n \}$ is a \defin{succession} of $\pi$ if
		$a-1$ and $a$ (or \(\overline{a-1}\) and \(\overline{a}\)) are in the same block. 
		If a partition has no succession, 
		then it is called a \defin{separated partition}.
	\end{definition}
	
	We use the notation $\defin{\MB(\pi)} \coloneqq \{m_i:  \pi_i \text{ is merging for } 2\le i\le k\}$,
	and $\defin{\Succeq(\pi)} \coloneqq \{a: \pm a \text{ is a succession}\}$.
	Further, 
	\[
	\defin{\mb(\pi)} \coloneqq |\MB(\pi)|, \quad
	\defin{\sueq(\pi)} \coloneqq |\Succeq(\pi)|, \quad\text{ and }\quad
	\defin{\nmb(\pi)} \coloneqq  k-\mb(\pi).
	\]
	\begin{example}
		The set partition $\pi=\typeBPartition{1, -10 | 2, 5 | 3, 4 | 6, -7, -8 | 9}$ 
		has $\MB(\pi)=\{6, 9\}$ and $\Succeq(\pi)=\{4, 8\}$.
	\end{example}

	In Table~\ref{tab:mbDistrib}, we present the first few values of the number of type $B$ set 
	partitions over $\langle n\rangle$ without zero block having $r$ merging blocks (or successions).
	\begin{table}[!ht]
		\centering
		\renewcommand{\arraystretch}{1.1}
		\begin{tabular}{c *{7}{r} r}
			\toprule
			$n \backslash r$ & \multicolumn{1}{c}{0} & \multicolumn{1}{c}{1} & \multicolumn{1}{c}{2} & \multicolumn{1}{c}{3} & \multicolumn{1}{c}{4} & \multicolumn{1}{c}{5} & \multicolumn{1}{c}{6} & $\qquad w_n$ \\
			\midrule
			1 &     1 &     &     &     &     &     &     &      1 \\
			2 &     2 &    1 &     &     &     &     &     &      3 \\
			3 &     6 &    4 &    1 &     &     &     &     &     11 \\
			4 &    24 &   18 &    6 &    1 &     &     &     &     49 \\
			5 &   116 &   96 &   36 &    8 &    1 &     &     &    257 \\
			6 &   648 &  580 &  240 &   60 &   10 &   1 &     &   1539 \\
			7 &  4088 & 3888 & 1740 &  480 &   90 &  12 &   1 &  10299 \\
			\bottomrule
		\end{tabular}
		\caption{The distribution of the number of merging blocks or successions over \(\typeBnzb{n}\).}\label{tab:mbDistrib}
	\end{table}
	
	\bigskip 
	\begin{example}\label{ex:sucMbDistr}
		
		The $11$ set partitions in $\typeBnzb{3}$ have the following joint distribution of $(\sueq,\mb)$.
		\begin{center}
			\begin{tabular}{@{}p{0.1\textwidth}p{0.1\textwidth}p{0.1\textwidth}p{0.1\textwidth}p{0.1\textwidth}p{0.1\textwidth}@{}}
				\toprule
				$(0,0)$ & $(1,0)$ & $(0,1)$ & $(1,1)$ & $(2,0)$ & $(0,2)$ \\
				\midrule
				\begin{tabular}[t]{@{}l@{}}
					\(\typeBPartition{1,-3 | 2} \)\\
					\(\typeBPartition{1,3 | 2}\) \\
					\(\typeBPartition{1,-2,3}\)
				\end{tabular}
				&
				\begin{tabular}[t]{@{}l@{}}
					\(\typeBPartition{1,-2,-3}\) \\
					\(\typeBPartition{1,2,-3}\)
				\end{tabular}
				&
				\begin{tabular}[t]{@{}l@{}}
					\(\typeBPartition{1 | 2,-3}\) \\
					\(\typeBPartition{1,-2 | 3}\)
				\end{tabular}
				&
				\begin{tabular}[t]{@{}l@{}}
					\(\typeBPartition{1 | 2,3}\) \\
					\(\typeBPartition{1,2 | 3}\)
				\end{tabular}
				&
				\begin{tabular}[t]{@{}l@{}}
					\(\typeBPartition{1,2,3}\)
				\end{tabular}
				&
				\begin{tabular}[t]{@{}l@{}}
					\(\typeBPartition{1 | 2 | 3}\)
				\end{tabular}
				\\
				\bottomrule
			\end{tabular}
		\end{center}
	\end{example}
	Let $\defin{\mathcal{U}_n^a} \coloneqq \{\pi \in \typeBnzb{n,k} : a\in\MB(\pi)\}$
	and $\defin{\mathcal{V}_n^a} \coloneqq \{\pi \in \typeBnzb{n,k}: a\in\Succeq(\pi)\}$.
	
	We saw in Example~\ref{ex:sucMbDistr} that $(\sueq,\mb)$ form a symmetric joint distribution.
	That is, for every pair $i,j$ the two sets
	\[
	\{ \pi \in \typeBnzb{n} : \sueq(\pi) = i \text{ and } \mb(\pi)=j \} \text{ and }
	\{ \pi \in \typeBnzb{n} : \sueq(\pi) = j \text{ and } \mb(\pi)=i \}
	\]
	have the same cardinality. This was proved in the type $A$ case in \cite{Beyeneetal},
	using a type of ``swap operation'' which we now define.
	
	\begin{definition}[Swap map]
		Let $\pi=\pi_1 \mid \pi_2 \mid \dotsb \mid \pi_k$ and let $i,j \in [k]$ and $a \in [n]$.
		\begin{equation}
			\defin{\Swap_a^{(i,j)}(\pi)}  \coloneqq \pi \text{ if $i=j$ or $\pm a \notin \pi_i\cup \pi_j$.}
		\end{equation}
		
		Otherwise, we let $I_a$ be the maximal integer interval in $\pi_i\cup \pi_j$ that starts with $\pm a$,
		and we move the elements of $I_a$ lying in $\pi_i$ to $\pi_j$ and vice versa.

		\medskip
		
		We now define the following maps.
		
		First, consider $\pi=\pi_1\mid \pi_2\mid \cdots \mid \pi_k\in \mathcal{U}_n^a$. 
		Let $a=\min(|\pi_i|)$, $1\le i\le k$. 
		Then $a{-}1~ (\text{or }\overline{a{-}1})\in \pi_j$ for some $j$.
		Note then that $\min(|\pi_j|)\le a-1<\min(|\pi_i|)$, whence $j<i$.
		
		Define the map $\mu_a: \mathcal{U}_n^a\mapsto \mathcal{V}_n^a$ by $\mu_a(\pi)=\pi'$, where $\pi'$ is obtained from $\pi$
		as follows.
		Let $\pi^*$ be the set partition obtained by merging the blocks $\pi_{i-1}$ and $\pi_i$ if $(a-1)\in \pi_j$
		or by merging the blocks $\pi_{i-1}$ and $-\pi_i$ if $(\overline{a-1})\in \pi_j$.
		Then put $\pi' = \Swap_a^{(i-1,j)}(\pi^*)$.
		We note that $\pm a$ becomes a succession of $\mu_a(\pi)$ and it is evident that $\NMB(\pi')=\NMB(\pi)$. 
	\end{definition}
	
	\medskip

	\begin{example}
		Let $\pi=\typeBPartition{1, 3, -5, 7, 10| 2, 4 | 6, -8| 9}$. 
		We have $\MB(\pi)=\{6, 9\}$, $\Succeq(\pi)=\emptyset$. 
		If $a=6$, then $i=3$, $j=1$ and $\pi^*=\typeBPartition{1, 3, -5, 7, 10| 2, 4, -6, 8 | 9}$. 
		Thus, $\pi'=\mu_6(\pi)=\typeBPartition{1, 3, -5, -6, 7, 10| 2, 4, 8 | 9}\in\mathcal{V}_{10}^6$. 
		Note that $\MB(\pi')=\{9\}, \Succeq(\pi')=\{6\}$, and that $\NMB(\pi')=\NMB(\pi)=\{1, 2\}$.
	\end{example}
	
	Next, consider $\pi=\pi_1\mid \pi_2\mid \cdots \mid \pi_k\in\mathcal{V}_n^a$. 
	Then $\pm(a-1, a)\in \pi_i$ for some $i$. 
	
	Define the map $\rho_a: \mathcal{V}_n^a\mapsto \mathcal{U}_n^a$ by $\rho_a(\pi)=\pi'$, where $\pi'$ is obtained from $\pi$ as follows. 
	Let $j$ be the smallest positive integer such that the 
	elements $\pm1, \pm2, \ldots, \pm (a-1)$ are in the first $j$ blocks of $\pi$. 
	Apply $\Swap_a^{(i,j)}$ to $\pi$ and then split the modified block $\pi_j$ before $\pm a$. 
	If $a < 0$, then replace the new block with its negative.
	\smallskip
	
	We note that if $\pm a$ is a succession, then $a$ becomes 
	the minimal element of a merging block of $\rho_a(\pi)$.
	Moreover, $\NMB(\pi')=\NMB(\pi)$.
	
	\begin{example}
		Let $\pi=\typeBPartition{1, -3, -4, 6, 9 | 2, -5, 8 | 7 | 10}$ with $\MB(\pi)=\{10\}$, $\Succeq(\pi)=\{4\}$, and let $a=4$. So $i=1$, $j=2$, and
		$\Swap_4^{(1,2)}(\pi)=\typeBPartition{1, -3, -5, 6, 9 | 2, -4, 8 | 7 | 10}$.
		Hence, we have $\typeBPartition{1, -3, -5, 6, 9 | 2 | -4, 8 | 7 | 10}$, and thus $\rho_4(\pi)=\pi'=\typeBPartition{1, -3, -5, 6, 9 | 2 | 4, -8 | 7 | 10}\in\mathcal{U}_{10}^4$.
		Note that $\MB(\pi')=\{4, 10\}, \Succeq(\pi')=\emptyset$ and that $\NMB(\pi')=\NMB(\pi)=\{1, 2, 7\}$.
	\end{example}
	
	For any $\pi\in\Pi_n^0$, $R =\{a_1, \dotsc, a_m\}\subseteq\MB(\pi)$ and $S=\{b_1, \ldots, b_s\}\subseteq\Succeq(\pi)$,
	we define $\defin{\psi_{R,S}(\pi)} \coloneqq \pi'$, where $\pi'$ is the set partition of type $B$ obtained
	from $\pi$ by applying $\mu_a$ for each element $a$ of $R$ and applying $\rho_b$ for each element $b$ of $S$.
	Thus,
	\[
	\psi_{R,S}=\mu_{a_1}\cdots\mu_{a_m}\rho_{b_1}\cdots\rho_{b_s}.
	\]
	It can easily be seen that $\psi_{R,S}$ is a bijection that exchanges
	the number of merging blocks and the number of successions.
	Therefore, we have the following theorem:
	\begin{theorem}
		For any $n \geq 1$, we have the following $qt$-symmetry:
		\[
		\sum_{\pi\in\typeBnzb{n}} q^{\mb(\pi)} t^{\sueq(\pi)}r^{\nmb(\pi)}
		=
		\sum_{\pi\in\typeBnzb{n}} t^{\mb(\pi)} q^{\sueq(\pi)}r^{\nmb(\pi)}.
		\]
		
	\end{theorem}
	\begin{corollary}
		The number of type $B$ set partitions without zero block over $\langle n\rangle$ having $r$ merging
		blocks and the number of type $B$ set partitions without zero block
		over $\langle n\rangle$ having $r$ successions are equal.
	\end{corollary}
	\begin{proposition}
		The number $w_k(n,s)$ of type $B$ set partitions without zero blocks over $\langle n\rangle$ having $k$ blocks and $s$ successions satisfies the recurrence relation
		\begin{equation}
			w_k(n,s) = w_{k-1}(n{-}1,s)+w_k(n{-}1,s{-}1)+(2k-1)w_k(n{-1},s),~ n\ge2,
		\end{equation}
		where $w_0(0,0)=1$, $w_1(0,0)=0$, $w_1(1,0)=1$, and $w_1(1,1)=0$.
		
		Further, the generating polynomials
		\begin{equation}
			Q_{n,s}(x)=\sum_{k\ge0}w_k(n,s)x^k
		\end{equation}
		satisfy
		\[
		Q_{n,s}(x)
		=(x-1)Q_{n-1,s}(x)+Q_{n-1,s-1}(x)+2xQ_{n-1,s}'(x).
		\]
		For \(0\leq s\leq n-1\), the polynomial \(Q_{n,s}(x)\) is
		real-rooted.
	\end{proposition}
	\begin{proof}
		Any type $B$ set partition over $\langle n\rangle$ having $k$ blocks and $s$ successions is
		obtained either from a type $B$ set partition over $\langle n-1\rangle$ having $k-1$ blocks or from a
		type $B$ set partition over $\langle n-1\rangle$ having $k$ blocks.
		
		In the former case, we insert $n$ in a singleton block. This implies that the number
		of blocks increases by one and the number of successions remains the same.
		This contributes $w_{k-1}(n{-1},s)$ to $w_k(n,s)$.
		
		In the latter case, we insert $\pm n$ in any of the existing $k$ blocks.
		If we insert $\pm n$ in the block containing $\pm (n{-}1)$, the number of
		successions increases by one. So this contributes $w_k(n-1,s-1)$ to $w_k(n,s)$.
		
		If we insert $\pm n$ in the block which does not contain $\pm (n{-}1)$, the number of
		successions remains the same. There are $2(k-1)$ ways to do this.
		Also, in the block containing $\pm (n{-}1)$, we can insert $\mp n$ without
		affecting the number of successions.
		Hence, this contributes $2(k-1)w_k(n-1,s)+w_k(n-1,s)=(2k-1)w_k(n-1,s)$ to $w_k(n,s)$.

		Multiplying the recurrence by \(x^k\) and summing over \(k\) gives
		the displayed recurrence for \(Q_{n,s}(x)\), with the convention
		\(Q_{n,-1}(x)=0\).

		We next reduce all values of \(s\) to the case \(s=0\). We claim that
		\begin{equation}\label{eq:Q-successions-reduction}
			Q_{n,s}(x)=\binom{n-1}{s}Q_{n-s,0}(x),
			\qquad 0\leq s\leq n-1.
		\end{equation}
		This follows by induction on \(n\). The case \(s=0\) is immediate. If
		\(1\leq s\leq n-2\), the induction hypothesis and the recurrence give
		\begin{align*}
			Q_{n,s}
			={}&\binom{n-2}{s}
			\left((x-1)Q_{n-s-1,0}+2xQ'_{n-s-1,0}\right)
			+\binom{n-2}{s-1}Q_{n-s,0}\\
			={}&\left(\binom{n-2}{s}+\binom{n-2}{s-1}\right)Q_{n-s,0},
		\end{align*}
		which proves the claim by Pascal's identity. For \(s=n-1\), the same
		recurrence gives \(Q_{n,n-1}=Q_{n-1,n-2}=\dotsb=Q_{1,0}=x\), so the
		claim also holds in the boundary case.

		It remains to prove that \(Q_{m,0}\) is real-rooted. Since every
		nonempty partition has at least one block, write
		\[
			Q_{m,0}(x)=xR_m(x).
		\]
		We have \(R_1(x)=1\), and the recurrence for \(Q_{m,0}\) becomes
		\begin{equation}\label{eq:R-zero-successions}
			R_m(x)=(x+1)R_{m-1}(x)+2xR'_{m-1}(x).
		\end{equation}
		The recurrence also shows that \(R_m\) has degree \(m-1\) and a
		positive leading coefficient. We prove by induction that its zeros are
		simple and negative. The cases \(R_1(x)=1\) and \(R_2(x)=x+1\) are
		immediate. Suppose that \(m\geq3\) and that the \(d=m-2\) zeros of
		\(R_{m-1}\) are
		\(\alpha_1<\dotsb<\alpha_d<0\). At such a zero,
		\[
			R_m(\alpha_i)=2\alpha_iR'_{m-1}(\alpha_i).
		\]
		These values alternate in sign. Moreover,
		\(R_m(0)=R_{m-1}(0)>0\), while the positive leading coefficient
		gives the opposite signs at \(-\infty\) and \(\alpha_1\).
		The intermediate value theorem therefore gives one zero of \(R_m\)
		in each of
		\[
			(-\infty,\alpha_1),\,
			(\alpha_1,\alpha_2),\dotsc,\,
			(\alpha_d,0).
		\]
		These are all the zeros of \(R_m\), so they are simple and negative.
		Thus, \(Q_{m,0}=xR_m\) is real-rooted, and
		\eqref{eq:Q-successions-reduction} completes the proof.
	\end{proof}
	

	\subsubsection{Another definition of a merging block}
	
	Let $\pi=\pi_1\mid \pi_2\mid \cdots\mid \pi_k$ be a type $B$ set
	partition without zero block over $\langle n\rangle$.
	For $2 \leq i \leq k$, we say that a block $\pi_i$ is merging if
	\[
	\max(\pi_{i-1})<\min(\pi_i).
	\]
	We call a block $\pi_i$ with this definition a \defin{normal merging block}.
	
	For instance, in $\pi=\typeBPartition{1 | 2, -6, -8 | 3, -4 | 5, 7}$, the second block is a merging
	block while the fourth block is a normal merging block.

	\begin{proposition}
		The distribution of the number of merging blocks and the number of normal merging
		blocks are the same in type $B$ set partitions without zero block.
	\end{proposition}
	\begin{proof}
		Let $\pi$ be a type $B$ set partition without zero block over $\langle n\rangle$ having $k$ blocks and $m$ merging blocks.
		Let $\{t_1, t_2, \dotsc, t_m\}$ be the set of indices of the merging blocks of $\pi$.
		For any $i\in[m]$, we define the map $\eta_{t_i}: \pi \mapsto \pi'$, where $\pi'$ is obtained from $\pi$
		by moving all barred integers in $\pi_{t_i}$ to the block $\pi_{t_i-1}$.
		
		We now show that the block $t_i$ of $\pi'$ is a normal merging block.
		Let $\pi'=\pi' \mid \pi_2' \mid  \cdots \mid \pi_k'$.
		All integers in $\pi_{t_i}'$ are unbarred, and since the integers moved to $\pi_{t_i-1}'$ are all barred,
		$\max(\pi_{t_i-1})=\max(\pi_{t_i-1}')$.
		Thus, $\max(\pi_{t_i-1}')\le\max(|\pi_{t_i-1}|)<\min(|\pi_{t_i}|)=\min(\pi_{t_i}')$.
		Therefore, $\pi_{t_i}'$ is a normal merging block.
		
		The procedure of $\eta_{t_i}$ is invertible; if $\pi_{t_i}'$ is a normal merging block, then move all barred integers of $\pi_{t_i-1}'$
		with absolute value larger than the minimum integer of $\pi_{t_i}'$ to $\pi_{t_i}'$.
		Therefore, $\eta_i^{-1}(\pi')=\pi$, and $\eta_i$ is indeed a bijection.
		
		Furthermore, the map $\eta \coloneqq \eta_{t_m}\circ\eta_{t_{m-1}}\circ\cdots\circ\eta_{t_1}$ is a bijection
		that exchanges the number of merging blocks and the number of normal merging blocks.
	\end{proof}
	
	\begin{example}
		Consider $\pi=\typeBPartition{1, -2 | 3, -4, -5 | 6, 8 | 7}$.
		Then $\pi$ has two merging blocks, the second and the third.
		So, $\pi'=\eta(\pi)=\eta_3\circ\eta_2(\pi)=\eta_3(\typeBPartition{1, -2, -4, -5 | 3 | 6, 8 | 7})=\typeBPartition{1, -2, -4, -5 | 3 | 6, 8 | 7}$.
		Note that $\pi'$ has two normal blocks, the second and the third.
	\end{example}

	
	\subsection{Type \texorpdfstring{$B$}{B} merging-free partitions}\label{sec:Merging-free}
	
	In this section, we study type $B$ merging-free partitions---partitions without merging blocks.
	
	Let \defin{$\mathcal{M}_n^B$} denote the set of type $B$ merging-free partitions without zero block over $\langle n\rangle$.
	
	We present a bijection between type $B$ set partitions over $\langle n\rangle$ and
	type $B$ merging-free partitions without zero block over $\langle n{+}1\rangle$.
	
	Define the map $f:\Pi_n^B\mapsto\mathcal{M}_{n+1}^B$ given by $f(\pi)=\pi'$,
	where $\pi'$ is obtained from $\pi$ as follows.
	\begin{enumerate}
		\item Let $\pi^*=\pi_0\mid \pi_1\mid \pi_2\mid \dotsb \mid \pi_k$ be the partition
		obtained from $\pi$ by removing all the singleton blocks.
		\item For $i$ from $k$ to $0$:  move the minimal element of $\pi_i$ to the block $i+1$.
		\item If $0\neq a\in\pi_0$, then move $a$ to the rightmost block whose minimal element less than $a$, and bar it.
		\item If $a$ is a singleton in $\pi$, then insert it in the  rightmost block of minimal element less than $a$.
		\item Increase every integer by $1$ and set the resulting partition $\pi'$.
	\end{enumerate}
	
	\begin{example}
		Consider $\pi = \typeBPartition{0, 3, 6 | 1, -7, 9 | 2 | 4 | 5, 8, -10 | 11}$.
		Then by removing the singleton blocks, we have $\pi^* = \typeBPartition{0, 3, 6 | 1, -7, 9 | 5, 8, -10}$.
		Move the minimal elements of the blocks of $\pi^*$ for $i$ from $2$ to $0$ to get $\typeBPartition{3, 6 | 0, -7, 9 | 1, 8, -10 | 5}$.
		Move $3$ and $6$ to their corresponding rightmost blocks and bar them. So we
		obtain $ \typeBPartition{0, -7, 9 | 1, -3, 8, -10 | 5, -6}$.
		Then insert the singletons to the corresponding rightmost blocks. We
		get $ \typeBPartition{0, -7, 9 | 1, 2, -3, 4, 8, -10 | 5, -6, 11}$.
		Finally, increase every integer by $1$, 
		obtaining $f(\pi)=\pi'=\typeBPartition{1, -8, 10 | 2, 3, -4, 5, 9, -11 | 6, -7, 12}$.
		
	\end{example}
	
	\begin{lemma}\label{lem:mergefree}
		The map $f$ is well-defined.
	\end{lemma}
	\begin{proof}
		First, we show that $f(\pi)=\pi'\in\mathcal{M}_{n+1}^B$.
		Let $a$ be the minimal element of the $j$-th block of $\pi'$, where $j>1$.
		Then $a{-}1$ is the minimal element of some non-singleton block of $\pi$.
		This implies that there is at least one integer in the block $j{-}1$ of $\pi'$ that in absolute value is greater than $a$.
		Assume, for a contradiction, that there is no such integer.
		That is, $a$ is greater than the absolute values of each element of the block $j{-}1$ of $\pi'$.
		So, $a{-}1$ cannot be the minimal element of a non-singleton block of $\pi$.
		This is a contradiction.
		Therefore, every minimal element (except $1$) of a block of $\pi'$ is less than the maximal element of the preceding block.
		Hence, $\pi'$ is merging-free.

		Next, consider $\pi=\rho$, where $\pi, \rho \in \Pi_n^B$.
		Then $f(\pi)=\pi'=\rho'=f(\rho)$.
		Therefore, $f$ is well-defined.
	\end{proof}
	
	\begin{theorem}\label{thm:bijectionf}
		The map $f:\Pi_n^B\mapsto\mathcal{M}_{n+1}^B$ defined by $f(\pi)=\pi'$ is a bijection.
	\end{theorem}
	\begin{proof}
		We show that $f$ is injective and surjective. 
		
		For injectivity, suppose $\pi, \rho\in\Pi_n^B$ and $f(\pi)=f(\rho)$. We show that $\pi=\rho$. 
		Suppose that they are not equal. Let the block $i$ be the leftmost block of the respective partitions 
		on which they disagree. This is because either the two blocks differ only by signs or there is an 
		element belongs to the $i$-th block of the first partition while it is not in the $i$-th block of the other partition. 
		
		In the former case, such block cannot be the zero block or a singleton block. 
		Since the procedure changes the signs of elements in the zero block, such integers would 
		have different signs in the image partitions. This is a contradiction. 
		
		In the latter case, let an element $a$ be in the $i$-th block of $\pi$ and $a$ is in the $j$-th block of $\rho, j>i$. 
		We only consider the possibility that $i=0$ and $a$ is barred in $\rho$. 
		In this case, $\overline{a}$ is not the minimal element of its block in $\rho$ and 
		the procedure shifts the minimal element of the block containing $\overline{a}$ to the right. 
		So there is a block with minimal element less than $a$ in the right side of the 
		block containing $\overline{a}$ in $f(\rho)$. However, it is not the case in $f(\pi)$ 
		since $\overline{a}$ would be in the rightmost block of minimal element less than $a$. 
		This is a contradiction. 
		
		Hence, $\pi=\rho$, and $f$ is injective.
		
		For surjectivity, let $\rho\in\mathcal{M}_{n+1}^B$. Let us construct a partition $\pi$ from $\rho$ as follows. 
		\begin{enumerate}
			\item Decrease the absolute value of every integer $1$.
			\item For each block, remove every non-minimal unbarred element that is less than the minimal element of the following block. All non-minimal unbarred elements of the last block should be removed.
			\item For each block, move every barred element whose absolute value is less than the minimal element of the following block to a new zero block and unbar it. All barred elements in the last block should also be unbarred and moved to the zero block.
			\item If the resulting partition has \(t+1\) blocks, then\\ for $i$ from 1 to $t$:\\
			\phantom{x}\hspace{3ex} move the minimal element of the \(i\)-th block to the block $i-1$.
			\item Include the removed integers in singleton blocks of the resulting partition.
		\end{enumerate}
		Set the resulting partition $\pi$. Then $\pi$ is a type $B$ set partition over $\langle n\rangle$, and we ensure that $f^{-1}(\rho)=\pi$ by construction.
		
		Therefore, $f$ is a bijection.
	\end{proof}
	\begin{example}
		Consider $\rho=\typeBPartition{1, 8, 10 | 2, 3, -4, 5, 9, 11 | 6, -7, 12}$.
		We decrease the absolute value of every integer by 1 
		to get $\typeBPartition{0, 7, 9 | 1, 2, -3, 4, 8, 10 | 5, -6, 11}$.
		Then 2, 4 and 11 are removed, so we have $\typeBPartition{0, 7, 9 | 1, -3, 8, 10 | 5, -6}$. 
		Now we move $\typeBPartition{-3}$ and $\typeBPartition{-6}$ to a new zero block and unbar them. 
		Thus we have  $\typeBPartition{3, 6 | 0, -7, 9 | 1, 8, -10 | 5}$.
		Next, for $i$ from $1$ to $3$, move the minimal element of the \(i\)-th block to the block $i-1$ 
		to obtain $\typeBPartition{0, 3, 6 | 1, -7, 9 | 5, 8, -10 }$.
		Finally, by including the removed elements in singleton blocks,
		we obtain $\typeBPartition{0, 3, 6 | 1, -7, 9 | 2 | 4 | 6, 8, -10 | 11}=\pi$.
	\end{example}
	\begin{corollary}
		The number of merging-free type $B$ set partitions without zero block 
		over $\langle n\rangle$ is the same as the number of type $B$ set partitions over $\langle n{-}1\rangle$. \qed
	\end{corollary}

	The number $a_{n, k}$ of type $A$ merging-free partitions over $[n]$ 
	having $k$ blocks satisfies the recurrence relation:
	\begin{equation}\label{eq1}
		a_{n, k}=k a_{n-1, k}+(n-2)a_{n-2, k-1},\quad n\geq2, \, k\geq1,
	\end{equation}
	where $a_{0, 0}=1$, $a_{1, 0}=0$, $a_{1, 1}=1$, see \cite{Be-Ma}.

	The following proposition gives a type $B$ analog of \eqref{eq1}.
	\begin{proposition}
		The number $b_{n,k}$ of type $B$ merging-free partitions over $\langle n\rangle$ 
		having $k$ blocks satisfies the recurrence relation
		\begin{equation}\label{bldist}
			b_{n,k}=(2k) b_{n-1,k} + 2(n-2)b_{n-2,k-1},\quad n\ge2
		\end{equation}
		where $b_{0, 0}=1$, $b_{1, 0}=0$, $b_{1, 1}=1$.
	\end{proposition}
	\begin{proof}
		We prove that $b_{n,k}=2^{n-k}a_{n,k}$, where $a_{n,k}$ is the number of type $A$ merging-free partitions 
		over $[n]$ having $k$ blocks. Observe that any type $B$ merging-free partition $\pi$ can be 
		obtained from a type $A$ merging-free partition $\rho$ over $[n]$ by adding signs to its 
		non-minimal elements. Since there are $k$ minimal elements, $\pi$ can have negative signs 
		on any subset of the remaining $n-k$ elements. 
		This contributes $2^{n-k}$, and hence, we have $b_{n,k}=2^{n-k}a_{n,k}$. 
		Then the recurrence can be obtained by a simple algebraic manipulation.
	\end{proof}
	
	We now consider a subset of type $B$ merging-free partitions over $\langle n\rangle$ consisting of a partition $\pi=\pi_1 \mid  \pi_2\mid \cdots\mid \pi_k$ such that if there is a barred element in $\pi_i, 1\le i\le k-1$, then it is the maximal element of $\pi_i$, and all elements of $\pi_k$ are unbarred. \begin{example}
		If $n=4$, then the elements of this subset are $\typeBPartition{1, 2, 3, 4}, \typeBPartition{1, 3 | 2, 4}, \typeBPartition{1, -3 | 2, 4}, \typeBPartition{1, 4 | 2, 3}$, $\typeBPartition{1, -4 | 2, 3}, \typeBPartition{1, 2, 4 | 3}$, $\typeBPartition{1, 2, -4 | 3}, \typeBPartition{1, 3, 4 | 2}$, and $\typeBPartition{1, 3, -4 | 2}$.
	\end{example}
	
	It can be seen that the number $\tilde{b}_{n,k}$ of type $B$ merging-free 
	partitions over $\langle n\rangle$ in this subset having $k$ blocks satisfies the following relation.
	\begin{equation}\label{bldist1}
		\tilde{b}_{n,k}=k \cdot \tilde{b}_{n-1,k}+2(n-2)\tilde{b}_{n-2,k-1},\quad n\ge2
	\end{equation}
	where $\tilde{b}_{0, 0}=1$, $\tilde{b}_{1, 0}=0$, $\tilde{b}_{1, 1}=1$.

	\subsubsection{Real-rootedness and interlacing roots}
	
	\begin{definition}[Interlacing and alternating, see \cite{Wagner1992}]
		
		Let $f$ and $g$ be polynomials with positive leading coefficients and with real
		roots $\{f_i\}$ and $\{g_i\}$, respectively.
		We say that \(f\) \defin{interlaces} \(g\) if \(deg(f) + 1 = deg(g) = d\)
		and
		\[
		g_1 \le f_1\le g_2 \le f_2 \le \dotsb \le f_{d-1}\le g_d
		\]
		Moreover, we say that \(f\) \defin{alternates left of} \(g\) if \(deg(f) = deg(g) = d\) and
		\[
		f_1\le g_1 \le f_2\le \ldots \le f_d \le g_d.
		\]
		
		We say that \(f\) \defin{interleaves} \(g\) if either \(f\) interlaces \(g\) or \(f\) alternates left of \(g\).
		We write this as \(f\ll g\).
		By convention, \(0 \ll 0, 0 \ll h\) and \(h \ll 0\), whenever \(h\) is a polynomial with a positive leading coefficient.
	\end{definition}

	\begin{proposition}\label{prop:interlacingRoots}
		Let $C_1$, $C_2$ be fixed positive real numbers. 
		Suppose the $c_{n, k}$ are defined via
		\[
		c_{n, k} = C_1 \cdot k  \cdot c_{n-1, k} + C_2 \cdot (n-2)c_{n-2, k-1},\quad n\geq2, \, k\geq 1,
		\]
		and $c_{0, 0}$, $c_{1, 0}$ and $c_{1, 1}$ are fixed initial conditions,
		and $c_{0,j}=c_{1,1+j}=0$ whenever $j\geq 1$.
		
		Let $P_n(t) \coloneqq \sum_{k=0}^n c_{n,k} t^k$.
		Then 
		\[
		P_n(t) =C_1 t \cdot P_{n-1}'(t)  +  C_2 \cdot (n-2) t \cdot P_{n-2}(t) 
		\]
		with the initial conditions $P_0(t)=c_{0,0}$, $P_1(t)=c_{1,1}t + c_{1,0}$.
		
		Moreover, $P_n(t)$ has only real zeros, and the roots of $P_k$ interlace the roots of $P_{k+1}$ for all $k \geq 1$.
	\end{proposition}
	\begin{proof}
		The fact that the recursion for the $c_{n,k}$ implies the 
		recursion for $P_n(t)$ is straightforward,
		so we shall only prove the last property using induction over $n$.
		
		Let $n \geq 3$ and suppose that $P_{n-2} \ll P_{n-1}$. 
		By Rolle's theorem, we also know that $P_{n-1}'(t) \ll P_{n-1}(t)$.
		Using some results from \cite{Wagner1992}, we can conclude that 
		\begin{align*}
			P_{n-1}(t) &\ll C_1 t \cdot P_{n-1}'(t)  \\
			P_{n-1}(t) &\ll C_2 (n-2) t \cdot  P_{n-2}(t).
		\end{align*}
		From this, it follows that
		\[
		P_{n-1}(t) \ll  C_1 t \cdot P_{n-1}'(t) +  C_2 (n-2) t \cdot P_{n-2}(t)  = P_{n}(t)
		\]
		and we are done.
	\end{proof}
	
	\begin{corollary}
		We may apply Proposition~\ref{prop:interlacingRoots}
		to the recursions in \eqref{eq1},  \eqref{bldist} and \eqref{bldist1}.
		Thus, the generating polynomials for $a_{n,k}$, $b_{n,k}$ and $\tilde{b}_{n,k}$, respectively, are real-rooted.
	\end{corollary}

		
	\subsection{Type \texorpdfstring{$B$}{B} separated set partitions}\label{sec:Separated}
	
	Recall that separated set partition is a partition in which all integers in the same block are non-consecutive.
	
	Let \defin{$\mathcal{S}_n^B$} denote the set of type $B$ separated set partitions without zero block over $\langle n \rangle$.
	
	We give a bijection between type $B$ set partitions over $\langle n\rangle$ and type $B$ 
	separated set partitions without zero block over $\langle n{+}1\rangle$.

	
	We define the map $g:\Pi_n^B\mapsto\mathcal{S}_{n+1}^B$ by $g(\pi)=\pi'$, 
	where $\pi'$ is obtained from $\pi$ as follows. Let $\pi=\pi_0 \mid \pi_1 \mid \cdots \mid \pi_k$. We apply the following procedure to $\pi$.
	\begin{enumerate}
		\item For $a$ from $2$ to $n$:\\
		\phantom{x}\hspace{3ex} if $a$ is a succession in the zero block, then bar it if and only if $a{-}1$ is unbarred;\\
		\phantom{x}\hspace{3ex} if $a$ (or $\overline{a}$) is a succession in block $\pi_i, i\ge1$, then move $a$ (or $\overline{a}$) to the zero block. 
		\item Increase every integer in absolute value by 1. 
	\end{enumerate}
	
	\begin{example}
		Let $\pi=\typeBPartition{0, 2, 3, 5 | 1, -7, -8 | 4, 6, 9, 10}$. Since 3 is a succession in the zero block and \(2\) is unbarred, we bar \(3\). Note that $\typeBPartition{-8}$ and $10$ are successions.
		So we move them to the zero block.
		Then we increase every integer in absolute value by 1. Therefore, $\pi'=\typeBPartition{1, 3, -4, 6, -9, 11 | 2, -8 | 5, 7, 10}\in\mathcal{S}_{11}^B$.
	\end{example}
	We define $\pi=g^{-1}(\pi')$ as follows. 
	Let $\pi'=\pi_1'\mid \pi_2'\mid \cdots\mid \pi_k'\in\mathcal{S}_{n+1}^B$.
	\begin{enumerate}
		\item Decrease every integer in absolute value by 1. 
		\item For $a$ from $2$ to $n$:\\
		\phantom{x}\hspace{3ex}if $a{-}1$ and $a$ are in the first block and either of them is barred, then unbar them;
		\phantom{x}\hspace{3ex}if $a$ (or $\overline{a}$) is in the first block while $a{-}1$ (or $\overline{a{-}1}$) is in another block, then\\
		\phantom{x}\hspace{3ex} move $a$ (or $\overline{a}$) to the block containing $a{-}1$ (or $\overline{a{-}1}$).
	\end{enumerate}
		
	We note that
	\begin{lemma}
		For any partition $\pi\in\mathcal{S}_n^B$, we have 
		\begin{equation}
			(g\circ g^{-1})(\pi)=\pi.
		\end{equation}
		\qed
	\end{lemma}
	Therefore, we have that
	\begin{proposition}
		The map $g$ is a bijection.
		
		\qed
	\end{proposition}
	
	\begin{corollary}
		The number of type $B$ set partitions over $\langle n\rangle$ having $k$ blocks is equal to the
		number of type $B$ separated set partitions over $\langle n+1\rangle$ having $k$ blocks.
		
		\qed
	\end{corollary}
	
	
	\subsection{Type \texorpdfstring{$B$}{B} merging-free separated set partitions}
	
	We let \defin{$\mathcal{MS}_n^B$} denote the set of type $B$ merging-free separated set partitions of over $\langle n\rangle$, and let $\defin{\mathcal{MS}_{n,k}^B\coloneqq} \{\pi\in\mathcal{MS}_{n,k}^B: \blocks(\pi)=k\}$, and $\defin{u_{n,k}\coloneqq} |\mathcal{MS}_{n,k}^B|$.
	
	\begin{theorem}\label{thmnombsuc}
		The number $u_{n,k}$ satisfies the recurrence
		\begin{equation}
			\label{rrmergfreesep}
			u_{n,k}=(2k-1)u_{n-1,k}+2(n-2)u_{n-2,k-1}, n\ge2, 1\le k\le \left\lceil\frac{n}{2}\right\rceil,
		\end{equation}
		where $u_{0,0}=1, u_{1,1}=1$.
		The generating polynomials $P_n(t)\coloneqq \sum_{k} u_{n,k} t^k$ satisfy the recursion 
		\begin{equation}
			P_n(t) = 2t \cdot  P'_{n-1}(t) - P_{n-1}(t) + 2t(n-2) P_{n-2}(t).
		\end{equation}
		Moreover, all the $P_n(t)$ are real-rooted.
	\end{theorem}
	\begin{proof}
		Any $\pi\in\mathcal{MS}_{n,k}^B$ can be obtained in either of the following two cases. 
		
		First, let $\pi\in\mathcal{MS}_{n-1,k}^B$. 
		By inserting $\pm n$ in any of the blocks except in the block containing $\pm(n{-}1)$, we obtain a 
		merging-free separated partition over $\langle n\rangle$ having $k$ blocks. 
		This contributes $(2k{-}1)u_{n-1,k}$ to $u_{n,k}$. The class of partitions obtained 
		in this case are characterized by the fact that the removal of $\pm n$ leaves a merging-free separated partition.
		
		Next, let $\pi\in\mathcal{MS}_{n-2,k-1}^B$. For each $a\in[n{-}2]$, increase every integer 
		greater than $a$ (in absolute value) by 1; let the $i$-th block be the rightmost block containing any of 
		the integers $\pm1, \pm2, \ldots, \pm a$; split this block into two blocks after the rightmost 
		integer of $\pm1, \pm2, \ldots, \pm a$; then insert $\pm n$ in the $i$-th block and insert $a{+}1$
		with the same sign as $a$ in the block $i{+}1$ of the resulting partition. 
		
		If $a{+}1$ is inserted with a bar, then take the maximal integer sequence starting with $\overline{a{+}1}$ in 
		its block and swap the signs of the elements in the interval. Then set the resulting partition $\pi^*$. 
		Observe that $\pi^*$ is a type $B$ set partition. The operation creates a new block, and 
		ensures that $\pi^*$ is merging-free since the maximum (in absolute value) of the $i$-th block 
		of $\pi^*$ is $n$ and it is larger than $a{+}1$, the minimum (in absolute value) of the newly created block. 
		If $a$ is in the $j$-th block of $\pi^*$, then set $\pi'=\Swap_{a+2}^{(j,i{+}1)}(\pi^*)$. 
		Observe that $\pi'=\pi^*$ if $I_{a+2}$ is empty, and that the resulting partition is succession-free. 
		The class of partitions obtained in this procedure 
		are characterized by that the removal of $\pm n$ results in a partition with merging blocks. 
		Since the procedure is reversible, the map $\pi\mapsto\pi'$ is indeed a bijection. 
		Thus, we have $2(n{-}2)u_{n-2,k-1}$ such partitions since $a$ has $n{-}2$ options, 
		and $n$ has $2$ options. 
		Hence, the recurrence.
		
		\medskip 
		
		We shall prove the slightly stronger statement that $P_{n-1} \interl P_{n}$.
		We make the induction hypothesis that $P_{n-2} \interl P_{n-1}$. 
		Now, by Rolle's theorem,
		\[
		P'_{n-1} \interl P_{n-1} \text{ so } P_{n-1} \interl 2tP'_{n-1}.
		\]
		The induction hypothesis implies that $P_{n-1} \interl 2t(n-2)P_{n-2}$ whenever $n \geq 3$.
		
		We can also show that $2t(n-2)P_{n-2} - P_{n-1}$ has a positive leading coefficient.
		Observe that the leading coefficient of $2t(n-2)P_{n-2} - P_{n-1}$ is
		\[
		\begin{cases}
			2(n-2)u_{n-2,\frac{n-1}{2}}, &\text{ if }  n \text{ is odd};\\ 
			2(n-2)u_{n-2,\frac{n-2}{2}}-u_{n-1,\frac{n}{2}},& \text{ if } n \text{ is even}.
		\end{cases}
		\]
		
		Since the first is clearly positive, we consider only the second case. Using the relation (\ref{rrmergfreesep}), we have
		
		\begin{align*}
			u_{n-2,\frac{n-2}{2}}=&\left(2\left(\frac{n-2}{2}\right)-1\right)u_{n-3,\frac{n-2}{2}}+2(n-4)u_{n-4,\frac{n-2}{2}-1}\\
			=&(n-3)u_{n-3,\frac{n-2}{2}}+2(n-4)u_{n-4,\frac{n-4}{2}}
		\end{align*}
		
		and
		
		\[
		u_{n-1,\frac{n}{2}}=\left(2\left(\frac{n}{2}\right)-1\right)u_{n-2,\frac{n}{2}}+2(n-3)u_{n-3,\left(\frac{n}{2}-1\right)}=0+2(n-3)u_{n-3,\frac{n-2}{2}},
		\]
		
		because $\frac{n}{2}>\lceil\frac{n-2}{2}\rceil$.
		Thus, the coefficient becomes
		\begin{align*}
			&2(n-2)\left((n-3)u_{n-3,\frac{n-2}{2}}+2(n-4)u_{n-4,\frac{n-4}{2}}\right)-2(n-3)u_{n-3,\frac{n-2}{2}}\\
			&=\left(2(n-2)(n-3)-2(n-3)\right)u_{n-3,\frac{n-2}{2}}+4(n-2)(n-4)u_{n-4,\frac{n-4}{2}}\\
			&= 2(n-3)^2u_{n-3,\frac{n-2}{2}}+4(n-2)(n-4)u_{n-4,\frac{n-4}{2}}\\
			&>0 \text{ whenever } n>3.
		\end{align*}
		%
		
		So by studying the graphs (see Figure~\ref{fig:interlacing}), we can conclude that
		\[
		P_{n-1} \interl 2t(n-2)P_{n-2} - P_{n-1}.
		\]
		We can then conclude that $P_{n-1}$ interlaces the sum:
		\[
		P_{n-1} \interl  2tP'_{n-1} + \left( 2t(n-2)P_{n-2} - P_{n-1} \right) = P_n,
		\]
		and we are done.
		
		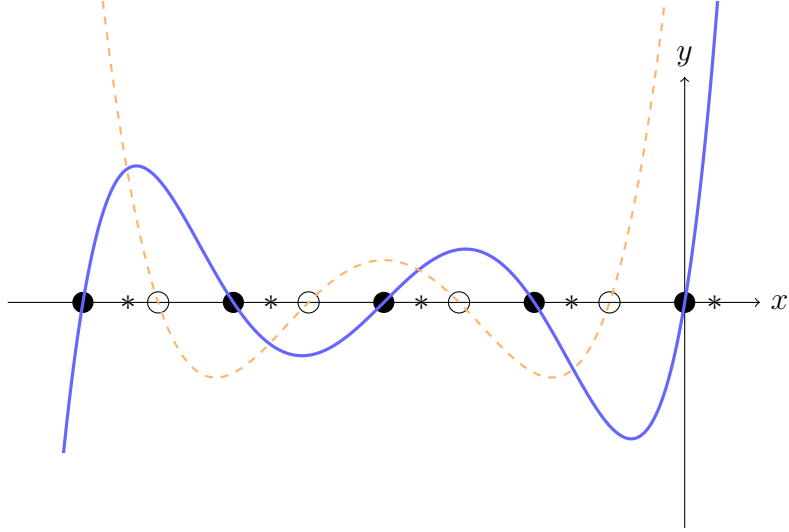
\begin{figure}[!ht]
			\begin{center}
				\begin{tikzpicture}[scale=2]
					\draw[->] (-4.5,0) -- (0.5,0) node[right] {$x$};
					\draw[->] (0,-1.5) -- (0,1.5) node[above] {$y$};
					
					\foreach \x in {0,-1,-2,-3,-4}
					\fill (\x,0) circle (2pt);
					\foreach \x in {-0.5,-1.5,-2.5,-3.5}
					\draw (\x,0) circle (2pt);
					
					\draw foreach \x in {0.2,-0.75,-1.75,-2.75,-3.7} {
						(\x,0) node {$\ast$}
					};
					
					\begin{scope}
						\clip (-4.5,-1) rectangle (0.5,2);
						
						\draw[blue!60, very thick, domain=-4.5:0.5, samples=200, smooth]
						plot (\x,{((\x)*(\x+1)*(\x+2)*(\x+3)*(\x+4))/4});
						
						\draw[orange!60, thick, dashed, domain=-4.5:0.5, samples=200, smooth]
						plot (\x,{((\x+0.5)*(\x+1.5)*(\x+2.5)*(\x+3.5))/2});
					\end{scope}
				\end{tikzpicture}
			\end{center}
			\caption{Illustration: The roots of $2t(n-2)P_{n-2}$ (marked $\bullet$) and $P_{n-1}$ (marked $\circ$).
				The roots of $P_{n-1}$ interlaces the roots of the difference $2t(n-2)P_{n-2} - P_{n-1}$
				(marked $\ast$).
			}\label{fig:interlacing}
		\end{figure}

		%
	\end{proof}

	\begin{example}
		The first few generating polynomials for $u_{n,k}$:
		\[
		\begin{array}{r}
			t \\
			t \\
			2 t^2+t \\
			10 t^2+t \\
			12 t^3+36 t^2+t \\
			140 t^3+116 t^2+t \\
			120 t^4+1060 t^3+358 t^2+t \\
			2520 t^4+6692 t^3+1086 t^2+t \\
		\end{array}
		\]
	\end{example}

	The number of type $B$ merging-free separated partitions  over $\langle n{+}1\rangle$ and the number of type $B$ set partitions without zero block over $\langle n\rangle$ are equal.
	
	Here, we present a bijection between these two sets.
	
	Define the map $h:\Pi_n^{\mathtt{0}} \mapsto \mathcal{MS}_{n+1}^B$ by $h(\pi)=\pi'$, where $\pi'$ is obtained from $\pi$ as follows. 
	\begin{enumerate}
		\item Let $\pi^*=\pi_1 \mid \pi_2\mid \cdots\mid \pi_k$ be the partition obtained from $\pi$ by removing the singleton blocks;
		\item For $i$ from $k$ to $1$:\\
		move $m_i=\min(|\pi_i|)$ to block $i+1$ (note that $m_k$ will form a new block);
		\item Increase every integer in absolute value by $1$, and then insert $1$ in the first block;
		\item For $a$ in increasing order:\\
		If $\pm a$ is a succession, then move $\pm a$ to the rightmost block of minimal element less than $a$;
		\item If $a$ is a singleton in $\pi$, 
		then insert $(a+1)$ in the rightmost block whose minimal element less than $(a+1)$, 
		ensuring that $a$ and $(a+1)$ have opposite signs.
	\end{enumerate}
	If $\pi$ and $\rho$ are two type \(B\) set partitions without zero block over $\langle n\rangle$ 
	and suppose that $\pi=\rho$, then one can see that $f(\pi)=f(\rho)$. 
	Thus, the map $h$ is well defined.
	
	We demonstrate the procedure of the map $h$ with an example.
	\begin{example}
		Consider $\pi=\typeBPartition{1, 2, 6 | 3 | 4 | 5, -9, -10 | 7, 8 | 11}\in\Pi_{11}^{\mathtt{0}}$. By removing the singletons, we obtain $\pi^*=\typeBPartition{1, 2, 6 | 5, -9, -10 | 7, 8}$. Then we apply the above procedure.
		\begin{enumerate}
			\item First, we shift the minimal elements and we obtain $\typeBPartition{ 2, 6 | 1, -9, -10 | 5, 8 | 7 }$.
			\item Next, we increase every integer in absolute value by 1, and then insert 1 in the first block to get $\typeBPartition{ 1, 3, 7 | 2, -10, -11 | 6, 9 | 8 }$.
			\item Next, we break any successions. There is only one succession $(\typeBPartition{ -11 })$,
			so we move $\typeBPartition{ -11 }$ to the rightmost block of minimal element less than $10$, the last block,
			and we get $\typeBPartition{1, 3, 7 | 2, -10 | 6, 9 | 8, -11}$.
			\item Finally, we insert the singletons. Since $3$ is not barred in the resulting partition,
			we insert $\typeBPartition{-4}$ in the rightmost block of minimal element less than $4$, the second block, and then insert $5$ also in the second block.
			Since $11$ is barred, then we insert $12$ in the last block.
			
			Therefore, $\pi'=\typeBPartition{1, 3, 7 | 2, -4, 5, -10 | 6, 9 | 8, -11, 12}\in\mathcal{MS}_{12}^B$.
		\end{enumerate}
	\end{example}
	\begin{lemma}\label{lemma:maph}
		For all $\pi\in\Pi_n^{\mathtt{0}}$, we have $h(\pi)\in\mathcal{MS}_{n+1}^B$.
		
		\qed
	\end{lemma}
	The procedure of the map $h$ is invertible. Let $\pi'\in\mathcal{MS}_{n+1}^B$ having $k$ blocks.
	We define $\pi=h^{-1}(\pi')$ as follows.
	\begin{enumerate}
		\item Delete 1 and decrease every integer in absolute value by 1.
		\item Consider the block \(i, 1\le i\le k\). For an element $\pm a$ different from the minimal element: 
		\begin{itemize}
			\item[-] if $a$ and $a{-}1$ have opposite signs, and $a$ is smaller than the minimal element of the block $i+1$, then remove $\pm a$ or;
			\item[-] if $a$ and $a{-}1$ have the same sign, and $a$ is smaller than the minimal element of the block $i{+}1$, then move $\pm a$ to the block containing $\pm(a{-}1)$;
		\end{itemize}
		\item For $i$ from $2$ to $k$:\\
		move the minimal element of block $i$ to block $i{-}1$;
		\item Insert the removed elements into the singleton blocks.
	\end{enumerate}
	\begin{example}
		Consider $\pi'=\typeBPartition{1, 3, 7 | 2, -4, 5, -10 | 6, 9 | 8, -11, 12}$. Deleting 1 and then decreasing every integer in absolute value by 1 gives us $\typeBPartition{2, 6 | 1, -3, 4, -9 | 5, 8 | 7, -10, 11}$. Since 3 and 2, 4 and 3, and, 11 and 10 have opposite signs, we remove 3, 4, and 11, and since 10 and 9 have the same sign, we move $\typeBPartition{-10}$ to the block containing $\typeBPartition{-9}$.
		From this we obtain $\typeBPartition{2, 6 | 1, -9,  -10 | 5, 8 | 7}$.
		Now, for \(i\) from 2 to 4: we move the minimal element of the \(i\)-th block to the block \(i{-1}\) to get $\typeBPartition{1, 2, 6 | 5, -9, -10 | 7, 8}$.
		Finally, we insert the removed elements into the singleton blocks.
		Thus, we have $\pi=\typeBPartition{1, 2, 6 | 3 | 4 | 5, -9, -10 | 7, 8 | 11}$.
	\end{example}
	\begin{lemma}\label{lemma:invertibleh}
		For any partition $\pi\in\Pi_n^{\mathtt{0}}$, we have 
		\begin{equation}
			(h^{-1}\circ h)(\pi)=\pi.
		\end{equation}
		
		\qed
	\end{lemma}
	
	By Lemma~\ref{lemma:maph} and Lemma~\ref{lemma:invertibleh}, we have the following theorem.
	\begin{theorem}
		The map $h$ is a bijection.
		\qed
	\end{theorem}
	
	\begin{proposition}
		Let $\pi\in\Pi_n^{\mathtt{0}}$ with $k$ non-singleton blocks
		If $h(\pi)=\pi'$, then $\pi'$ has $k+1$ blocks.
		\qed
	\end{proposition}
	
	\begin{corollary}\label{propnonsing}
		The number of type $B$ merging-free separated partitions 
		over $\langle n\rangle$ having $k$ blocks and the number of type $B$ set 
		partitions without zero block over $\langle n{-}1\rangle$ 
		having exactly $k{-}1$ non-singleton blocks are equal.
		\qed
	\end{corollary}
	
	Up to a shift on both $n$ and $k$, the number of type $A$ merging-free separated set partitions over $[n]$ having 
	$k$ blocks is the same sequence as OEIS entry number \cite[\seqnum{A008299}]{Sl},
	counting set partitions without singletons.
	In \cite{Beyeneetal}, the authors raised a question of a natural bijection between these sets.
	Therefore, we answer this question in the following corollary.
	\begin{corollary}
		The restriction of the bijection $h$ is a bijection between type $A$ set partitions 
		without singletons over $[n]$ and type $A$ merging-free separated set partitions over $[n+1]$.\qed 
	\end{corollary}
	
	\section{Flattened Stirling permutations}\label{sec:flatStirlingperm}
	Let \(\defin{[n]_2\coloneqq} \{1, 1, 2, 2, \dotsc, n, n\}\) denote the multiset with each element in \([n]\) appearing twice.
	We recall that a \defin{Stirling permutation} of order \(n\) is a word \(w = w_1 w_2\cdots w_{2n}\) on the multiset \([n]_2\),
	in which every letter in \([n]\) appears exactly twice and any values between the repeated instances of \(i\) must be greater than \(i\).
	The number of Stirling permutations of order \(n\) is given by the double factorial \((2n-1)!!\).
	
	The maximal contiguous weakly increasing subword of \(w\) is called a \defin{run} of \(w\). If the sequence of leading terms of the runs appear in weakly increasing order, then we say \(w\) is a \defin{flattened Stirling permutation} of order \(n\), see \cite{Bucketal}.
	
	Buck et al.~presented a bijection between flattened Stirling permutations of order $n{+}1$ and
	type \(B\) set partitions over \(\langle n\rangle\) \cite[Theorem 1.1]{Bucketal}.
	Since, by Theorem~\ref{thm:bijectionf}, type \(B\) set partitions over \(\langle n\rangle\) and
	type \(B\) merging-free set partitions over \(\langle n{+}1\rangle\) are in bijection, one can
	ask for a direct bijection between type \(B\) merging-free set partitions over \(\langle n\rangle\) and
	flattened Stirling permutations of order $n$.
	We will present one as follows.
	
	Let $\pi=\pi_1\mid \pi_2\mid \cdots\mid \pi_k$ be a type $B$ merging-free set partition over $\langle n\rangle$. Let $m_i=\min(|\pi_i|)$ for each $i$. We write the elements of each block $\pi_i$ that are larger in absolute values than $m_{i+1}$ such that the barred elements precede the unbarred ones. We construct a word $w=w_1w_2\cdots w_{2n}$ from $\pi$ as follows.
	\begin{enumerate}
		\item Replace 1 with $1~1$.
		\item If $a$ is in the $i$-th block with $a<m_{i+1}, i<k$, or is in the last block, and it is barred, then insert it between the 1's. All such \(a\) inserted between the 1s are arranged in increasing order of their absolute values.
		\item For each $i\ge1$, insert $m_{i+1}$ in the $i$-th block after the barred elements.
		\item Replace each element $a$ different from the minimal elements with $aa$.
		\item Remove the bars and the separators.
	\end{enumerate}
	\begin{example}
		Let $\pi=\typeBPartition{1, -6, 7 | 2, 3, -4, 8, -9, -10 | 5, -11, 12}$. Then we rewrite the elements of the blocks as $\pi=\typeBPartition{1, -6, 7 | 2, 3, -4, -9, -10, 8 | 5, -11, 12}$. Then the procedure is as follows.
		\begin{enumerate}
			\item $\typeBPartition{1, 1, -6, 7 | 2, 3, -4, -9, -10, 8 | 5, -11, 12}$
			\item $\typeBPartition{1, -4, -11, 1, -6, 7 | 2, 3, -9, -10, 8 | 5, 12}$
			\item $\typeBPartition{1, -4, -11, 1, -6, 2, 7 | 2, 3, -9, -10, 5, 8 | 5, 12}$
			\item $\typeBPartition{1, -4, -4, -11, -11, 1, -6, -6, 2, 7, 7 | 2, 3, 3, -9, -9, -10, -10, 5, 8, 8 | 5, 12, 12}$
			\item $w=\typeBPartition{1, 4, 4, 11, 11, 1, 6, 6, 2, 7, 7, 2, 3, 3, 9, 9, 10, 10, 5, 8, 8, 5, 12, 12}$
		\end{enumerate}
	\end{example}
	\begin{lemma}
		Let \(\pi\) be a type \(B\) merging-free set partition over \(\langle n\rangle\), and \(w\)
		be the word obtained from \(\pi\) by the above procedure. Then \(w\) is a flattened Stirling permutation of order \(n\).
	\end{lemma}
	\begin{proof}
		We show that \(w\) satisfies the required conditions for being a flattened Stirling permutation of order \(n\). 
		
		Each \(a\in[n]\) appears twice in \(w\), and that the minimal element of a block of \(\pi\) begins a run in \(w\). Only the instances of these minimal elements may have values between them. Clearly, any values appearing between the two instances of 1 (if any) are larger in absolute value than 1. 
		
		Let \(m_i=\min(|\pi_i|), i>1\), denoting the minimal element of the \(i\)-th block of \(\pi\).
		If all elements of \(\pi_{i-1}\) which are in absolute value larger than \(m_i\) are barred, then we have the instances \(m_im_i\) with no values between them. This satisfies the condition for a Stirling permutation.
		
		If, on the other hand, \(\pi_{i-1}\) contains at least one unbarred element larger than \(m_i\), then \(m_i\) is inserted after the barred elements of \(\pi_{i-1}\) whose absolute values exceed \(m_i\); that is, before unbarred elements larger than \(m_i\). Consequently, the values appearing between the two instances of \(m_i\) are larger than \(m_i\), which again satisfies the defining condition of a Stirling permutation.
		
		Therefore, \(w\) is a Stirling permutation. 
		
		Moreover, since the sequence of the leading terms of the runs of \(w\) consists of the minimal elements of the blocks of \(\pi\), the sequence is weakly increasing. Hence, \(w\) is a flattened Stirling permutation. 
	\end{proof}
	\begin{theorem}
		The correspondence \(\pi\mapsto w\), where \(w\) is the word obtained from \(\pi\) by the above procedure, is a bijection.
	\end{theorem}
	\begin{proof}
		Since the sets of type \(B\) merging-free set partitions over \(\langle n\rangle\) and flattened Stirling permutations of order $n$ have the same cardinality, it suffices to show that the correspondence \(\pi\mapsto w\) is one-to-one. 
		
		Let \(\pi\) and \(\beta\) be type \(B\) merging-free set partitions over \(\langle n\rangle\), and \(\pi\mapsto w\) and \(\beta\mapsto w'\). Assume that \(\pi\neq\beta\). Then there are different cases to consider.
		
		First, let \(\pi\) and \(\beta\) differ only in sign. That is, there is an element \(a\) that is barred in \(\pi\) and unbarred in \(\beta\), or vice versa. Note that \(\pm a\) is not the minimal element of a block of both partitions. In the partition where \(a\) is barred (say \(\pi\)), if the block following the block containing \(\overline{a}\) has a minimal element larger than \(a\), or \(\overline{a}\) lies in the last block, then the procedure moves \(\overline{a}\) to the position between the 1's. So, in this case, since the procedure does not move \(a\) in \(\beta\), the two image words are clearly different. If the block following the block containing \(\overline{a}\) has a minimal element \(m\) smaller than \(a\), then the procedure inserts \(m\) after the barred elements of the block containing \(\overline{a}\). Since \(a\) is barred, \(aa\) would be to the left of the appearances of \(m\) in \(w\). However, since \(a\) is unbarred in \(\beta\), and \(a>m\), we have \(aa\) between the appearances of \(m\) in \(w'\). Therefore, \(w\neq w'\).
		
		Now, consider the smallest \(i\) such that an element \(a\) belongs to the \(i\)-th block of, say \(\pi\), but it does not belong to the \(i\)-th block of \(\beta\). We only consider the case where \(a\) is barred in both partitions, and the procedure moves \(\overline{a}\) to the position between the 1's. Let \(\overline{a}\) belongs to the \(j\)-th block of \(\beta\). 
		
		If \(\overline{a}\) belongs to the last block of the partitions, then because of the minimality of \(i\) and the assumption that \(\pi\neq\beta\), \(\beta\) has necessarily more number of blocks than \(\pi\). Thus, we must have \(w\neq w'\). 
		
		Assume that \(\overline{a}\) is not in the last block of \(\pi\). For each \(i\), let \(m_i\) and \(m_i'\) be the minimal elements of the \(i\)-th blocks of \(\pi\) and \(\beta\), respectively. Then we have that \(m_j'< a< m_{i+1}\). Since \(\pi\) and \(\beta\) agree on the first \(i-1\) blocks, \(\pm m_j'\) and \(\overline{a}\) are in the same block of \(\pi\). If \(m_j'\) is barred in \(\pi\), then the procedure moves it to the position between the instances of 1's while it is not the case in \(\beta\). If \(m_j'\) is unbarred in \(\pi\), then all integers between the instances of \(m_i\) and \(m_j'\) will be smaller than \(m_j'\) (if there are any) and we have the instance \(m_j'm_j'\) with no integer between them in \(w\). However, since \(\beta\) is merging-free, we must have at least one integer larger than \(m_j'\) between the instances of \(m_i'\) and \(m_j'\) or between the instances of \(m_j'\) and \(m_j'\) in \(w'\). As a result, we have \(w\neq w'\).
		
		Therefore, the correspondence \(\pi\mapsto w\) is one-to-one, and hence, it is indeed a bijection.
	\end{proof}
	
	\section{Flattened signed permutations}
	\label{sec:flatsignedperm}
	Let $\pi=\pi_1\mid \pi_2\mid \cdots\mid \pi_k$ be any type $B$ set partition over $\langle n\rangle$.
	Then we define  $\sigma\coloneqq Flatten(\pi)$, a signed permutation 
	over $\defin{[\pm n]} \coloneqq \{\pm1, \pm2, \ldots, \pm n\}$
	obtained from $\pi$ by removing the block separators.
	For instance, if $\pi=\typeBPartition{1, -8 | 2, 4, 7 | 3, 6, -9 | 5}$, then we have $\sigma=Flatten(\pi)=\typeBPartition{1, -8, 2, 4, 7, 3, 6, -9, 5}$.
	
	Let $\defin{\mathcal{R}_n^B}$ denote the set of signed permutations obtained by 
	flattening all type $B$ merging-free set partitions without zero block 
	over $\langle n\rangle$. Let $\defin{d_n}$ be the cardinality of $\mathcal{R}_n^B$. 
	Note that this number is the shifted Dowling number.
	
	We give the definition of a modular run for a signed permutation as follows.
	\begin{definition}
		Let $\sigma=\sigma_1\sigma_2\cdots\sigma_n\in \mathcal{R}_n^B$, then we define a \defin{modular run} (or in short \defin{mrun})
		of $\sigma$ as a maximal sequence $\sigma_i\sigma_{i+1}\cdots\sigma_{i+j}$ such
		that $|\sigma_i|<|\sigma_{i+1}|<\cdots<|\sigma_{i+j}|$.
		If $\sigma_i\sigma_{i+1}\cdots\sigma_{i+j}$ is an mrun of $\sigma$, then we
		say that $\sigma_i$ is its \defin{bottom} and $\sigma_{i+j}$ its \defin{top}.
	\end{definition}
	
	For instance, $\typeBPartition{1, -8}, \typeBPartition{2, -4, 7}, \typeBPartition{3, 6, -9}$ and $5$ are
	the mruns of $\sigma=\typeBPartition{1, -8, 2, -4, 7, 3, 6, -9, 5}$.
	\begin{remark}\label{rmsignperm}
		Note that a signed permutation is in $\mathcal{R}_n^B$ if and only
		if the absolute value of the top of each mrun is larger than the bottom of the following mrun.
	\end{remark}
	Given a signed permutation $\sigma=\sigma_1\sigma_2\cdots\sigma_n$, we say that $i\in[n{-}1]$ 
	is a \defin{descent} if $\sigma_i>\sigma_{i+1}$ and an \defin{ascent} otherwise.
	The \defin{major index} of $\sigma$ is given by $\maj(\sigma)=\sum_{i\in\Des(\sigma)}i$, 
	where $\defin{\Des(\sigma)}$ is the set of descents of $\sigma$. 
	We let \(\defin{\des(\sigma)}\coloneqq |\Des(\sigma)|\).

	We study the distribution of descents and major index on the class of signed permutations $\mathcal{R}_n^B$.
	
	The following table presents the first few values of $d_{n,k}\coloneqq |\{\sigma\in\mathcal{R}_n^B: ~\des(\sigma)=k\}|$.
	\begin{table}[H]
		\centering
		$\begin{array}{l|llllllll}
			n\backslash k&0&1&2&3&4&5&6&7\\
			\hline
			1&1 \\
			2&1&1 \\
			3&1&4&1 \\
			4&1&11&11&1 \\
			5&1&26&62&26&1 \\
			6&1&57&266&266&57&1 \\
			7&1&120&991&1864&991&120&1\\
			8&1&247&3405&10667&10667&3405&247&1\\
		\end{array}$
		\caption{The values of $d_{n,k}$, for $1\le n\le8, 0\le k\le n-1$}
	\end{table}
	
	\begin{proposition}
		\label{propsymmetry}
		The distribution of descents over the set $\mathcal{R}_n^B$ is symmetric, that is,
		\begin{equation}
			d_{n,k}=d_{n,n-k-1}.
		\end{equation}
	\end{proposition}
	\begin{proof}
		Let $\sigma\in\mathcal{R}_n^B$, and let $\sigma'$ be obtained from $\sigma$ by exchanging the signs
		of all the integers except for the bottoms of the mruns of $\sigma$, clearly $\sigma'\in\mathcal{R}_n^B$.
		Then $\sigma$ has a descent at the position $i$ if and only if $\sigma'$ has an ascent in $i$.
		This is obvious if $\sigma_i$ is not the top of an mrun, because the signs of both $\sigma_i$ and $\sigma_{i+1}$ are changed in $\sigma'$.
		
		Suppose then that $\sigma_i$ is the top of an mrun (except the last one), hence $\sigma_{i+1}$ is the bottom of the following mrun and therefore, $\sigma_{i+1}'=\sigma_{i+1}>0$. If $\sigma_i<0$, then $\sigma$ has an ascent at $i$ while $\sigma'$ has a descent at $i$, because of Remark \ref{rmsignperm} regarding permutations of $\mathcal{R}_n^B$, whereas if $\sigma_i>0$, then $\sigma$ has a descent at $i$ (for the same reason) while $\sigma'$ has an ascent at $i$. Therefore, $\sigma$ has $k$ descents if and only if $\sigma'$ has $n-k-1$ descents.
		
	\end{proof}
	For instance, if $\sigma=\typeBPartition{1, 3, -7, 2, 6, 4, 5, -9, 8}$, then $\sigma'=\typeBPartition{1, -3, 7, 2, -6, 4, -5, 9, 8}$. Then $\sigma$ has 3 descents and $\sigma'$ has $5=9-3-1$ descents.
	\begin{corollary}
		The distribution of the major index over the set $\mathcal{R}_n^B$ is symmetric.
	\end{corollary}
	\begin{proof}
		If $\sigma\in\mathcal{R}_n^B$ and $\sigma'$ is obtained from $\sigma$ by exchanging the signs of all the integers except for the bottoms of the mruns of $\sigma$, then $\maj(\sigma)=m$ if and only if $\maj(\sigma')=\binom{n}{2}-m$.
	\end{proof}
	\begin{corollary}
		Let $n\ge1$ and
		\[
		G_n(q,t)=\sum_{\sigma\in\mathcal{R}_n^B}q^{\des(\sigma)}t^{\maj(\sigma)}.
		\]
		Then
		\begin{equation}
			G_n(q,t)=q^{n-1}t^{\binom{n}{2}}G_n\left(q^{-1},t^{-1}\right).
		\end{equation}\qed
	\end{corollary}

	\begin{theorem}
		The number $d_{n,k}$ of permutations in $\mathcal{R}_n^B$ having $k$ descents satisfies the identity:
		\begin{equation}\label{desrunsort}
			d_{n,k}=\sum_{m=1}^{\min\{k{+}1,\lceil \frac{n}{2}\rceil\}}2^{m-1}\binom{n-2m+1}{k-m+1}a_{n,m},
		\end{equation}
		where $a_{n,m}$ is the number of run-sorted permutations over $[n]$ having $m$ runs.
	\end{theorem}
	\begin{proof}
		Any permutation of $\mathcal{R}_n^B$ having $k$ descents can be obtained from a run-sorted permutation $\tau=\tau_1\tau_2\cdots\tau_n$ over $[n]$ having $m$ runs, where $m\le k+1$, by adding some bars to the appropriate integers. Notice that $\tau$ is also in $\mathcal{R}_n^B$ and has precisely $m-1$ descents.
		
		When one adds a bar to either $\tau_n$ or to an integer that is neither the minimum nor the maximum of a run of $\tau$, then the number of descents increases by 1. Therefore, when we put bars on $k-m+1$ of such integers, we obtain a permutation with $k$ descents. There are $\binom{n-2m+1}{k-m+1}$ such choices. 
		
		On the other hand, when one adds a bar to the maximum of a run of $\tau$ (except the last), then the number of descents remains the same. Therefore, one can choose freely any subset of such maxima, add bars to them and still obtain a permutation with $k$ descents. There are $2^{m-1}$ such choices.
		
		Therefore, by the product rule, and then taking the sum over all possible $m$, we have the result.
	\end{proof}
	\begin{corollary}
		Let $b_{n,m}$ be the number of type $B$ merging-free partitions over \(\langle n\rangle\) without zero block having $m$ blocks. Then
		\begin{equation}
			d_{n,k}=\sum_{m=1}^{\min\{k+1,\lceil \frac{n}{2}\rceil\}}2^{-(n-2m+1)}\binom{n-2m+1}{k-m+1}b_{n,m},
		\end{equation}
	\end{corollary}
	\begin{proof}
		This follows from \eqref{desrunsort} using the relation 
		\begin{equation}
			\label{abrelation}
			b_{n,m}=2^{n-m}a_{n,m}.
		\end{equation}
	\end{proof}
	\begin{theorem}
		For all integers $n$ and $k$, $0\le k\le n-1$, we have
		\begin{equation}
			\label{rrdes}
			d_{n,k}=d_{n-1,k-1}+d_{n-1,k}+2\sum_{i=2}^{n-1}\binom{n-2}{i-1}\sum_{j=0}^{i-2}\binom{i-2}{j}d_{n-i,k-j-1}.
		\end{equation}
	\end{theorem}
	\begin{proof}
		
		Let $\sigma=\sigma_1\sigma_2\cdots\sigma_n$ be a flattened signed permutation having $k$ descents. Any flattened signed permutation $\sigma$ has $\sigma_1=1$ and there are different possibilities for $\sigma_2$. 
		
		\textbf{Case-1:} Let \(\sigma_2=\pm2\). If $\sigma_2=2$, then let $\sigma'$ be a permutation obtained from $\sigma$ by deleting $\sigma_1=1$ and then reducing every integer in absolute value by 1. Then $\sigma'\in\mathcal{R}_{n-1}^B$ and $\des(\sigma')=\des(\sigma)$. If $\sigma_2=\overline{2}$, then let $\sigma'$ be a permutation obtained from $\sigma$ by deleting $\sigma_1=1$, unbarring $2$, and then reducing every integer in absolute value by 1. Then $\sigma'\in\mathcal{R}_{n-1}^B$ and $\des(\sigma')=\des(\sigma)-1$. These explain the first two terms of the right-hand side of \ref{rrdes}.
		
		\textbf{Case-2:} Let $\sigma_2\neq\pm 2$. Then there exists an index $i>2$ such that $\sigma_{i+1}=2$ in $\sigma$. By definition, the prefix $1\sigma_2\cdots\sigma_{i}$ in absolute value is an increasing sequence of $i$ integers, and $2$ is always unbarred. If we delete this prefix and standardize the remaining part, we obtain a permutation $\tau\in\mathcal{R}_{n-i}^B$ with descents preserved except for those in the prefix. 
		Let us examine the number of descents we would have reduced while deleting the prefix $1\sigma_2\cdots\sigma_{i}$. The number of bars in the word $\sigma_2\cdots\sigma_{i-1}$ equals the number of descents in $1\sigma_2\cdots\sigma_{i-1}$. Also, $\sigma_{i}$ forms a descent with $\sigma_{i+1}=2$ if $\sigma_{i}$ is unbarred or it becomes a foot of a descent if $\sigma_{i}$ is barred.
		
		For the reverse of this procedure, we construct the prefix for each $\tau\in\mathcal{R}_{n-i}^B$ as follows. We choose $i-1$ integers from $n-2$ integers in $\binom{n-2}{i-1}$ ways because \(1\) and \(2\) cannot be chosen. Then we choose $j$ integers from the $i-2$ (excluding the largest) integers in $\binom{i-2}{j}$ ways to bar them. Hence we have an mrun consisting of 1 and those integers, where \(j\) of them are barred. Let \(1\sigma_2\sigma_3\cdots\sigma_i\) be this mrun. Note that it has $j$ descents. For each $\tau\in\mathcal{R}_{n-i}^B$, there is a unique permutation \(\sigma\in\mathcal{R}_{n}^B\), where \(1\sigma_2\sigma_3\cdots\sigma_i\) is the first mrun of \(\sigma\). Therefore, if $\tau$ has $k-j-1$ descents, then this contributes $\sum_{j=0}^{i-2}\binom{i-2}{j}d_{n-i,k-j-1}$. Together with the descent created by barring or unbarring $\sigma_{i}$, we have $k-j-1+j+1=k$ descents. 
		This explains the last part of the right-hand side of \eqref{rrdes}.
	\end{proof}

	\subsection{Valley-hopping}
	We discuss several permutation statistics which we will use in the rest of this paper.
	
	Given a signed permutation $\sigma=\sigma_1\sigma_2\cdots\sigma_n$, we say that $\sigma_i$ is
	\begin{itemize}
		\item a \defin{peak} if $\sigma_{i-1}<\sigma_i>\sigma_{i+1}$;
		\item a \defin{valley} if $\sigma_{i-1}>\sigma_i<\sigma_{i+1}$;
		\item a \defin{double ascent} if $\sigma_{i-1}<\sigma_i<\sigma_{i+1}$, where ($2\le i\le n-1$);
		\item a \defin{double descent} if $\sigma_{i-1}>\sigma_i>\sigma_{i+1}$, where ($2\le i\le n-1$).
	\end{itemize}
	Define $\Pk(\sigma), \Val(\sigma), \Dasc(\sigma)$, and $\Ddesc(\sigma)$ to be the set of peaks, valleys, double ascents, and double descents of $\sigma$, respectively. Moreover, let $\pk(\sigma)\coloneqq|\Pk(\sigma)|, \val(\sigma)\coloneqq|\Val(\sigma)|, \dasc(\sigma)\coloneqq|\Dasc(\sigma)|$, and $\ddesc(\sigma)\coloneqq|\Ddesc(\sigma)|$.

	Here, we are using the conventions $\sigma_0=-\infty, \sigma_{n+1}=\infty$. We note that $\Pk(\sigma), \Val(\sigma)\subseteq\{\sigma_1, \sigma_2, \ldots, \sigma_n\}$ and $\Dasc(\sigma), \Ddesc(\sigma)\subseteq\{\sigma_2, \sigma_3, \ldots, \sigma_{n-1}\}$.
	\begin{example}
		Let $\sigma=\typeBPartition{1, 3, 8, 2, -6, -7, 4, -5}$.
		Then $\Pk(\sigma)=\{4, 8\}, \Val(\sigma)=\{\overline{5}, \overline{7}\}, \Dasc(\sigma)=\{3\}$, and $\Ddesc(\sigma)=\{2, \overline{6}\}$.
	\end{example}
	We shall soon define a group action on $\mathcal{R}_n^B$ induced by involutions. 
	First recall  a $\mathbb{Z}_2^n$-action on the symmetric group $\mathfrak{S}_n$
	which is commonly known as the \defin{modified Foata--Strehl action} or \defin{valley-hopping}. 
	This action is based on a classical group action of Foata and Strehl \cite{Fo-St}, 
	was independently introduced by Shapiro, Woan, and Getu~\cite{Shapiroetal}, 
	and later rediscovered by Bränd\'{e}n~\cite{Branden2004}. An interested reader can also consider a paper by N. Lafreni\'ere and Y. Zhaung \cite{La-Zh}. 
	Fix a permutation $\sigma\in\mathfrak{S}_n$ and a letter $x\in[n]$. 
	We may write $\sigma=w_1w_2xw_4w_5$, where $w_2$ is the maximal consecutive subword immediately to the left
	of $x$ whose letters are all smaller than $x$, and $w_4$ is the maximal consecutive subword 
	immediately to the right of $x$ whose letters are all smaller than $x$; this decomposition is
	called the \defin{$x$-factorization of $\sigma$}.
	\medskip 

	Define $\varphi_x:\mathfrak{S}_n\mapsto\mathfrak{S}_n$ by
	\[
	\defin{\varphi_x(\sigma)} \coloneqq\begin{cases}
		w_1w_4xw_2w_5, & \text{ if $x$ is a double ascent or double descent of } \sigma,\\
		\sigma, & \text{ if $x$ is a peak or valley of } \sigma.
	\end{cases}
	\]
	
	Here, the conventions $\sigma_0=\sigma_{n+1}=\infty$ are used. 
	It is easy to see that $\varphi_x$ is an involution, and that for all $x, y\in[n]$, 
	the involutions $\varphi_x$ and $\varphi_y$ commute with each other. 
	Given a subset $S\subseteq[n]$, we define $\varphi_S:\mathfrak{S}_n\rightarrow\mathfrak{S}_n$ 
	by $\defin{\varphi_S} \coloneqq\Pi_{x\in S} \varphi_x$. 
	The involutions $\{\varphi_S\}_{S\subseteq[n]}$ define a $\mathbb{Z}_2^n$-action on $\mathfrak{S}_n$.
	
	We are now ready to define the group action on $\mathcal{R}_n^B$. 
	Fix a permutation $\sigma=\sigma_1\sigma_2\cdots\sigma_n\in\mathcal{R}_n^B$. 
	Note that a letter $x$ of $\sigma$ is either the bottom of an mrun or not, and that $x$ is always unbarred 
	if it is the bottom of an mrun. 
	\medskip 
	
	Define $\psi_i:\mathcal{R}_n^B\mapsto\mathcal{R}_n^B$ by 
	\[
	\defin{\psi_i(\sigma)} \coloneqq\begin{cases}
		\sigma_1\sigma_2\cdots\overline{\sigma_i}\cdots\sigma_n, & \text{ if $\sigma_i$ is 
		not the bottom of an mrun of } \sigma,\\
		\sigma, & \text{ if $\sigma_i$ is the bottom of an mrun of } \sigma.
	\end{cases}
	\]
	It is easy to see that $\psi_i$ is an involution, and that for all $i, j\in[n]$, 
	the involutions $\psi_i$ and $\psi_j$ commute with each other. 
	Given a subset $S\subseteq[n]$, we define $\psi_S:\mathcal{R}_n^B\rightarrow\mathcal{R}_n^B$ by 
	$\defin{\psi_S} \coloneqq\Pi_{i\in S}\psi_i$. 
	For instance, if $\sigma=\typeBPartition{1, 3, -8, 2, -6, 7, 4, 5}$ and $S=\{2,3,7\}$, 
	then $\psi_S=\typeBPartition{1, -3, 8, 2, -6, 7, 4, 5}$. 
	The involutions $\{\psi_S\}_{S\subseteq[n]}$ define a $\mathbb{Z}_2^n$-action on $\mathcal{R}_n^B$.

	\begin{proposition}
		For $S=[n]$ and $\sigma\in\mathcal{R}_n^B$, we have
		\begin{enumerate}
			\item $\dasc(\psi_S(\sigma))=\ddesc(\sigma)$;
			\item $\ddesc(\psi_S(\sigma))=\dasc(\sigma)$.
		\end{enumerate}
	\end{proposition}
	\begin{proof}
		The proof follows from the proof of Proposition \ref{propsymmetry}.
	\end{proof}
	We say that a subset $\Pi\subseteq\mathcal{R}_n^B$ is \defin{invariant} 
	under valley-hopping if for every $S\subseteq[n]$ and a 
	permutation $\sigma\in\mathcal{R}_n^B$, we have $\psi_S(\sigma)\in\Pi$ 
	(equivalently, if $\Pi$ is a disjoint union of orbits of the valley-hopping action).

	\subsection{Gamma-positivity}
	
	Let $\mathcal{R}_{n,k}^B$ be the set of permutations in $\mathcal{R}_n^B$ having $k$ mruns.
	Also, given $\tau\in\mathcal{R}_n^B$, we let $\defin{\Orb(\tau)} \coloneqq \{\psi_S(\tau): S\subseteq[n]\}$
	denote the orbit of $\tau$ under valley-hopping.
	
	\begin{proposition}\label{prop:numbOrbits}
		The number of orbits of $\mathcal{R}_n^B$ under the valley-hopping
		is given by the Bell number $B_{n-1}$.
	\end{proposition}
	\begin{proof}
		The number of orbits under a group action is the average number of fixed elements.
		Since valley-hopping defines a $\mathbb{Z}_2^n$-action on $\mathcal{R}_n^B$, it has cardinality $2^n$,
		equal to the number of subsets of $[n]$. Let $S\subseteq[n]$ and $\Fix(\psi_S)$ 
		denote the set of all fixed elements of $\mathcal{R}_n^B$ under $\psi_S$, i.e.,~$\defin{\Fix(\psi_S)}\coloneqq \{\sigma\in\mathcal{R}_n^B : \psi_S(\sigma)=\sigma\}$.
		So we prove that
		\[
		\frac{1}{2^n}\sum_{S\subseteq[n]}|\Fix(\psi_S)|=B_{n-1}.
		\]
		We note that $\psi_S$ fixes $\sigma$ if and only if $\sigma_i$ is the bottom of an mrun of $\sigma$ for all $i\in S$.
		
		Let $S=\{i_1, i_2, \ldots, i_{k-1}\}$, where $i_1\neq1$.
		If $T=\{1\}\cup S$, then $|\Fix(\psi_T)|=|\Fix(\psi_S)|$. Also, we  note that $S'\subseteq S$ implies that $|\Fix(\psi_{S'})|=|\Fix(\psi_S)|$. We have $|\Fix(\psi_S)|=b_{n,k}$, where $b_{n,k}$ is the number of permutations in $\mathcal{R}_n^B$ having $k$ mruns, if $S$ is a subset of non-consecutive integers of the set $\{3, 4, \ldots, n\}$. Otherwise, $|\Fix(\psi_S)|=0$.
		Also, any such $S$ contributes the same number of fixed elements associated to any of its non-empty subset.
		So there are $2^{k-1}-1$ fixed elements associated to the non-empty subsets of $S$. Since we need to consider $S$ or $T$, they contribute 2 to the number of fixed elements.
		
		Moreover, we note that $|\Fix(\psi_{\emptyset})|=|\Fix(\psi_{\{1\}})|=|\mathcal{R}_n^B|=d_n$.
		
		Recall \eqref{abrelation} that $b_{n,k}=2^{n-k}a_{n,k}$, where $a_{n,k}$ is the number of run-sorted
		permutations over $[n]$ having $k$ runs, and that $\sum_{k=1}^{\lceil\frac{n}{2}\rceil}a_{n,k}=B_{n-1}$.
		
		Therefore, we have 
		\begin{align*}
			\sum_{S\subseteq[n]}|\Fix(\psi_S)|=&2\left(d_n+\sum_{k=2}^{\lceil\frac{n}{2}\rceil}(2^{k-1}-1)b_{n,k}\right)\\
			=&2\left(d_n+\sum_{k=2}^{\lceil\frac{n}{2}\rceil}2^{k-1}b_{n,k}-\sum_{k=2}^{\lceil\frac{n}{2}\rceil}b_{n,k}\right)\\
			=&2\left(d_n+2^{n-1}\sum_{k=2}^{\lceil\frac{n}{2}\rceil}a_{n,k}-\left(d_n-2^{n-1}\right)\right)\\
			=&2\left(2^{n-1}(B_{n-1}-1)+2^{n-1}\right)\\
			=&2^nB_{n-1}
		\end{align*}
		as desired.
	\end{proof}
	
	\begin{lemma}\label{lemmaorbit}
		Let $\tau\in\mathcal{R}_{n,k}^B$.
		Then
		\begin{equation}
			\sum_{\sigma\in\Orb(\tau)}x^{\des(\sigma)}= (2x)^{k-1}(1+x)^{n-2k+1}.
		\end{equation}
	\end{lemma}
	\begin{proof}
		For every $\sigma\in\Orb(\tau)$, we have $k-1\le\des(\sigma)\le n-k$.
		The minimum number of descents is attained when $\sigma$ is a permutation such that if there is a barred element,
		then it is only the top element of one of its mruns (except the last one).
		This explains the factor $2^{k-1}$. Other additional descents are created by adding bars to any of the eligible $n-2k+1$ elements.
	\end{proof}
	In light of Proposition~\ref{prop:numbOrbits}, the number of orbits under
	valley-hopping having $k$ mruns is equal to $a_{n,k}$, the number
	of run-sorted permutations over $[n]$ having $k$ runs.
	Thus, the following lemma follows from Lemma~\ref{lemmaorbit}.
	
	\begin{lemma}
		We have
		\begin{equation}\label{eqn:lemmaorbit1}
			\sum_{\sigma\in\mathcal{R}_{n,k}^B}x^{\des(\sigma)}=	a_{n,k}(2x)^{k-1}(1+x)^{n-2k+1}.
		\end{equation}
		\qed
	\end{lemma}
	
	Let $\tilde{b}_{n,k}$ denote the number of type $B$ merging-free
	partitions with $k$ blocks as before in \eqref{bldist1}.
	Since we have $\tilde{b}_{n,k}=2^{k-1}a_{n,k}$, we then have the
	following theorem by summing \eqref{eqn:lemmaorbit1} over all orbits.
	\begin{theorem}\label{thm:gammapos}
		For all $n\ge 0,$ 
		\begin{equation}
			\sum_{\sigma\in\mathcal{R}_{n+1}^B} x^{\des(\sigma)}
			=\sum_{i=0}^{\lceil n/2 \rceil} \tilde{b}_{n+1,i+1} \cdot x^i (1+x)^{n-2i}.
		\end{equation}
		\qed
	\end{theorem}
	Note that this is an instance of so-called \defin{\(\gamma\)-positivity}.
	For an introduction to this topic, see the survey by Athanasiadis~\cite{Athanasiadis2018}.

	\begin{corollary}
		Since $\sum_k \tilde{b}_{n,k} t^k$ has only (negative) real roots, 
		it follows from \cite{Branden2004} that the polynomials 
		\[
		\sum_{\sigma\in\mathcal{R}_{n+1}^B} x^{\des(\sigma)}
		\]
		are also real-rooted.
	\end{corollary}

	\subsection{\texorpdfstring{\(\des\)}{des} is homomesic under valley-hopping}
	
	A statistic is \defin{$r$-mesic} under a group action if its average value over each orbit is equal to $r$ (see \cite{La-Zh}).
	A statistic is \defin{homomesic} if it is $r$-mesic for some $r$.
	
	We show that the statistic \(\des\) is \((n-1)/2\)-mesic.
	In other words, the average number of descents in each valley-hopping orbit
	is \((n-1)/2\).
	
	\begin{theorem}
		The descent statistic is $(n-1)/2$-mesic under valley-hopping.
	\end{theorem}
	\begin{proof}
		Let $\tau\in\mathcal{R}_{n,k}^B$. Then we have $|\Orb(\tau)|=2^{n-k}$.
		Let
		\[
			F_\tau(x)\coloneqq
			\sum_{\sigma\in\Orb(\tau)}x^{\des(\sigma)}
			=(2x)^{k-1}(1+x)^{n-2k+1}.
		\]
		The average number of descents in this orbit is
		\[
			\frac{F_\tau'(1)}{F_\tau(1)}
			=(k-1)+\frac{n-2k+1}{2}
			=\frac{n-1}{2},
		\]
		as desired.
	\end{proof}
\subsection*{Acknowledgments.} 
The second author acknowledges the financial support provided by the Department of Informatics 
at Universit\'e Paris Cit\'e and IRIF. We appreciate the hospitality we got from IRIF laboratory 
during the research visit of the second author.

\printbibliography
	
\end{document}